\documentclass[a4paper,12pt]{article}

 \usepackage{amsmath, amssymb, amsthm}

\usepackage{color}
\usepackage{graphicx}
\usepackage{graphics}

\definecolor{vio}{rgb}{0.5,0,0.5}
\definecolor{gre}{rgb}{0.1,0.6,0}
\definecolor{ora}{rgb}{0.8,0.2,0.1}

\numberwithin{equation}{section}

\newtheorem{theorem}{Theorem}[section]
\newtheorem{proposition}[theorem]{Proposition}

\newtheorem{definition}[theorem]{Definition}

\newtheorem{hypothesis}[theorem]{Hypothesis}

\newfont{\ctv}{msam10}

\newcommand{\beq}{\begin{equation}}
\newcommand{\eeq}{\end{equation}}

\newcommand{\bbox}{\mbox{\ctv \symbol{4}}}
\def\QED{{$\hfill\bbox$}}
\newenvironment{pf}[1]{\par\vskip1mm{\noindent\it #1.}\ }{\QED\par\vskip2mm}

\def\real{\mathbb{R}}
\def\nat{\mathbb{N}}

\def\ve{\varepsilon}
\def\vr{\varrho}
\def\vp{\varphi}
\def\la{\lambda}
\def\La{\Lambda}

\def\thre{\theta_{\rm ref}}

\def\expe{\mathrm{e}}

\def\for{\mbox{ for }\ }
\def\ale{\mbox{ a.e.}}
\def\dd{\,{\rm d}}

\def\on{^{(n)}}

\def\bw{\mathbf{w}}
\def\bfv{\mathbf{v}}
\def\bfD{\mathbf{D}}
\def\bve{\boldsymbol{\ve}}

\def\FF{\mathcal{F}}
\def\UU{\mathcal{U}}
\def\SS{\mathcal{S}}
\def\KK{\mathcal{K}}
\def\DD{\mathcal{D}^*}

\def\ipi{\int_0^{1}\!}

\def\stop{\mathfrak{s}}
\def\play{\mathfrak{p}}

\def\sumk{\sum_{k=1}^n}
\def\sumkj{\sum_{k=1}^{n-1}}

\newcommand{\limn}{\lim\limits_{n\to \infty}}

\def\be{\begin{equation}\label}
\def\ee{\end{equation}}

\def\beal{\begin{align}}
\def\endal{\end{align}}

\def\barr{\begin{array}}
\def\earr{\end{array}}

\def\bearr{\begin{eqnarray}}
\def\eearr{\end{eqnarray}}

\def\bears{\begin{eqnarray*}}
\def\eears{\end{eqnarray*}}

\def\bpf{\begin{pf}}
\def\epf{\end{pf}}

\def\no{\nonumber}

\usepackage{wasysym}

\begin{document}

\title{{A new phase field model for material fatigue in oscillating elastoplastic beam}
\footnote{Supported by GA\v CR Grant P201/10/2315, and by the institutional support
for the development of research organizations I\v C 47813059 and 67985840. The work of Michela Eleuteri
was supported by the FP7-IDEAS-ERC-StG
Grant \#256872 (EntroPhase)}}
\author{
Michela Eleuteri\footnote{Dipartimento di Matematica ``F. Enriques'', Universit\`a degli Studi di Milano,
via Saldini 50, 20133 Milano, Italy. Dipartimento di Matematica e Informatica ``U. Dini'', Universit\`a degli Studi di Firenze, viale Morgagni 67/a, 50134 Firenze, Italy. 
E-mail {\tt Michela.Eleuteri@unimi.it, eleuteri@math.unifi.it}.},
Jana Kopfov\'a\footnote{Mathematical Institute of the Silesian University, Na
Rybn\' i\v cku 1, 746 01 Opava, Czech Republic, E-mail {\tt
Jana.Kopfova@math.slu.cz}.},
and Pavel Krej\v{c}\'{\i}
\footnote{Institute of Mathematics,
Academy of Sciences of the Czech Republic, \v{Z}itn\'a 25, CZ-11567 Praha 1, Czech Republic,
E-mail {\tt krejci@math.cas.cz}.}
}

\maketitle

\begin{abstract}
We pursue the study of fatigue accumulation in an oscillating elastoplastic beam
under the additional hypothesis that the material can partially recover
by the effect of melting. 
The full system consists of the momentum and energy balance equations, an evolution
equation for the fatigue rate, and a differential inclusion for the phase dynamics.
The main result consists in proving the existence and uniqueness of a strong solution.
\end{abstract}

\section*{Introduction}\label{intr}

{It was shown in \cite{ksb} that the Kirchhoff-Love method of reducing the 3D problem
of transversal oscillations of a solid elastoplastic beam with the single yield von Mises plasticity law
leads to the beam equation with a multiyield hysteresis Prandtl-Ishlinskii constitutive operator.
The present authors have used in \cite{EKK13-2} {(see also \cite{EKK12}, \cite{KS})} the Prandtl-Ishlinskii formalism to propose a model for the cyclic fatigue accumulation in an oscillating beam and to study its properties; { results have been  obtained correspondingly also for the plate, see \cite{EKK13-1}, \cite{EKK14-1}, \cite{EKK14-2}.}
Here, we extend the model by taking into account the possibility of partial fatigue recovery
by the effect of melting when a solid-liquid phase transition takes place.}

{The fatigue accumulation law is still based on the observation that there exists a proportionality
between accumulated fatigue and dissipated energy, see \cite{bdk,fla}. Unlike in \cite{EKK13-2}
and similarly as in \cite{fla}, we assume that out of all dissipative components in the energy
balance, only the purely plastic dissipation produces damage. This makes the mathematical problem easier: {the system} of equations then does not develop singularities in finite time
and a unique regular solution is proved to exist on every bounded time interval. {On the other hand we consider here an additional difficulty  - we assume that  the weight function $\varphi$ in the definition of the Prandtl-Ishlinskii operator depends also on the fatigue parameter $m$; this has been considered also in \cite{EKK14-1} and \cite{EKK14-2}.
}

{The unknowns of the problem are the transversal displacement $w \in \real$ of the beam,
the absolute temperature $\theta > 0$, the fatigue variable $m \ge 0$, and the phase
variable $\chi\in [0,1]$. The full system of equations consists of the momentum balance equation
(the beam equation with a fatigue dependent hysteresis operator), the energy balance
equation, the fatigue accumulation equation and the phase transition equation.
The model is derived in detail in Section \ref{phys}.
}

{The problem is rigorously stated in Section \ref{stat}, where we also check the
thermodynamic consistency of the system and collect some preliminary material in Section 3. In Section 
\ref{form} we carry out formally the a priori estimates that allow us to construct the solution
of the full system. In Section \ref{appr}, we apply these ideas to a spatially discrete
scheme and derive estimates that are sufficient for proving that the space discrete
approximations converge to a solution of the original problem in appropriate function spaces.
The main existence and uniqueness Theorem \ref{t2} is proved in Section \ref{proo}.
}


\section{The model}\label{phys}

\subsection{Governing equations}\label{cons}

We consider a transversally inhomogeneous beam of length $1$, and denote by $x \in [0,1]$
the {longitudinal} variable, by $t \in [0,T]$
the time variable, by $w(x,t)$ the transversal displacement of the point $x$ at time $t$,
by $\ve(x,t) = w_{xx}(x,t)$ the linearized curvature, and by $\sigma(x,t)$ the bending moment.
We assume a thermo-visco-elasto-plastic scalar constitutive law in the form
\be{b3}
\sigma = B\ve + P[m,\ve] + \nu \ve_t - \beta(\theta - \thre)\,,
\ee
where {$B > 0$} is a constant hardening modulus, {$m\ge 0$ is a scalar time and
space dependent parameter describing the accumulation of fatigue, where $m=0$
corresponds to zero fatigue,} $P[m,\ve]$ is a fatigue dependent Prandtl-Ishlinskii
constitutive operator of elastoplasticity defined below in Subsection \ref{pran},
$\nu$ is the viscosity coefficient, $\beta$ is the thermal bending coefficient
related to a layered structure of the beam, $\theta>0$ is the absolute temperature,
and $\thre$ is a fixed referential temperature {(more specifically, the melting temperature)}.
Following \cite{ksb}, Newton's law of motion is formally written as
\be{motion}
\rho w_{tt} - \alpha w_{xxtt} +\sigma_{xx} = F(x,t),
\ee
where $\alpha = \rho l^2/12$ and $l>0$ is the thickness of the beam, $\rho$ the mass density and $F$
is the external load.

With the constitutive law \eqref{b3}, we associate the {\em free energy operator\/}
\be{b4}
\FF(\ve, \theta, \chi) = c\theta(1 - \log(\theta/\thre)) + \frac{B}{2} \ve^2 + V[m,\ve]
- \beta(\theta - \thre)\ve-\frac{L}{\thre}(\theta-\thre)\chi+I_{[0,1]}(\chi)\,,
\ee
where $V[m,\ve]$ is the fatigue dependent Prandtl-Ishlinskii potential \eqref{pot}, $c$
(the specific heat capacity) and $L$ (the latent heat) are given constants, and $I_{[0,1]}$
is the indicator function of the interval $[0,1]$. The {\em entropy operator\/}
$\SS$ and {\em internal energy operator\/} $\UU$ then read
\begin{align}\label{b5}
\SS(\ve, \theta, \chi) &= - \frac{\partial \FF}{\partial \theta} = c \log(\theta/\thre) + \beta\ve {+ \frac{L}{\thre}\chi}\,,
\\ \label{b6}
\UU(\ve, \theta, \chi) &= \FF(\ve, \theta)+ \theta \SS(\ve, \theta) = c\theta + \frac{B}{2} \ve^2 + V[m,\ve]
+ \beta\thre\ve + L\chi + I_{[0,1]}(\chi)\,.
\end{align}
We consider the first and the second principles of thermodynamics in the form
\begin{align}\label{b7}
\UU(\ve, \theta, \chi)_t +  q_x &= \ve_t \sigma + g\,,\\ \label{b8}
\SS(\ve, \theta, \chi)_t + \Big(\frac{q}{\theta}\Big)_x &\ge \frac{g}{\theta}\,,
\end{align}
where $q = - \kappa \theta_x$ is the heat flux with a constant heat conductivity $\kappa>0$,
and $g$ is the heat source density.
Note that \eqref{b7} is the energy conservation law, \eqref{b8} is the Clausius-Duhem inequality.

The evolution of the phase variable $\chi$ is governed by the inclusion
$-\gamma\chi_t \in \partial_\chi \FF$, that is,
\be{phase}
-\gamma\chi_t \in \partial I_{[0,1]}(\chi)-\frac{L}{\thre}(\theta-\thre),
\ee
where $\gamma > 0$ is a characteristic time of phase transition, and $\partial I_{[0,1]}$
is the subdifferential of the indicator function $I_{[0,1]}$. Indeed, we necessarily
have $\chi \in [0,1]$, and we interpret $\chi = 0$ as the solid phase,
$\chi = 1$ as liquid, and the intermediate values correspond to the relative liquid 
content in a mixture of the two.

Let $D[m,\ve]$ be the dissipation operator defined in \eqref{diss} associated
with the Prandtl-Ishlinskii operator $P[m,\ve]$. 
The analysis of the so-called {\em rainflow method of cyclic fatigue accumulation\/}
in elastoplastic materials carried out in \cite{bdk} has shown a close relation
between accumulated fatigue and dissipated energy, similarly as in \cite{fla}.
Here, we assume in addition that partial recovery of the damaged material is possible
under strong local melting. Mathematically, this is expressed in terms of the evolution equation
for the fatigue variable $m$
{
\be{fati}
m_t(x,t) \in - \partial I_{[0,\infty)}(m) -h(\chi_t(t)) + \int_0^1 \la(x-y) D[m,\ve](y,t)\dd y,
\ee
}
where $h$ is a nonnegative nondecreasing function, $\la$ is a nonnegative smooth function with (small) compact support and $ D[m,\ve]$ is the fatigue dependent dissipation operator, \eqref{diss}.
{The subdifferential $\partial I_{[0,\infty)}$ of the indicator function $I_{[0,\infty)}$
ensures that the fatigue parameter remains nonnegative.}

The meaning of \eqref{fati} is {simple. If no phase transition takes place or
if the material solidifies, that is, $\chi_t \le 0$, then fatigue
at a point $x$ increases proportionally to the energy dissipated in a
neighborhood of the point $x$. On the other hand, under strong
melting if $\chi$ grows faster than the plastic dissipation rate, the fatigue
may decrease until it possibly reaches the unperturbed state $m=0$.}


\subsection{{Hysteresis operators}}\label{pran}

Let us first recall the definition of the stop.

\begin{definition}\label{dpi1}
{
Let $u \in W^{1,1}(0,T)$ and a closed connected set $Z \subset \real$ be given.
The variational inequality
\be{stop}
\left.
\barr{ll}
u(t) = z(t)+\xi(t)& \forall t \in [0,T],\\[2mm]
z(t) \in Z & \forall t \in [0,T],\\[2mm]
\dot\xi(t)(z(t) - y) \ge 0 \ \ale & \forall y \in Z\,,\\[2mm]
z(0) = z^0 \in Z,
\earr
\right\}
\ee
defines the stop and play operators $\stop_Z$ and $\play_Z$ by the formula
\be{pi2s}
z(t) = \stop_Z[z^0,u](t)\,, \quad \xi(t) = \play_Z[z^0,u](t)\,.
\ee
}
For a canonical choice {of $Z= [-r,r]$ with some $r>0$ and for the initial condition
$z(0) = Q_r(u(0))$,} where $Q_r$ is the projection of $\real$ onto the interval $[-r,r]$, 
we simply write
\be{pi2}
z(t) = \stop_r[u](t)\,, \quad \xi(t) = \play_r[u](t)\,.
\ee
\end{definition}

A simple proof of the following easy properties of the play and stop can be found e.g. in \cite{K}.

{
\begin{proposition}\label{3p2}
Let $u_1, u_2 \in W^{1,1}(0,T)$, a closed connected set $Z \subset \real$, and data $z_1^0, z_2^0
\in Z$ be given, $z_i = \stop_Z[z_i^0, u_i]$, $\xi_i = u_i - z_i$, $i=1,2$.
Then
\begin{itemize}
\item[{\rm (i)}] $\displaystyle{(z_1(t) - z_2(t))
(\dot u_1(t) - \dot u_2(t))
\ \ge\ \frac12 \frac{\dd}{\dd t} (z_1(t) - z_2(t))^2}$
\quad a.\,e.;
\item[{\rm (ii)}] $\displaystyle{|\dot \xi_1(t) - \dot \xi_2(t)|
 + \frac{\dd}{\dd t}|z_1(t) - z_2(t)|
 \ \le\ |\dot u_1(t) - \dot u_2(t)|}$ \quad a.e.;
\item[{\rm (iii)}] $\displaystyle{|z_1(t) - z_2(t)|
\ \le\ |z_1^0 - z_2^0| + 2\,\max_{0\le\tau \le t}|u_1(\tau) - u_2(\tau)|}$ \quad
$\forall t \in [0,T]$;
\item[{\rm (iv)}] $\dot\xi_i(t)\dot u_i(t) = \dot\xi_i(t)^2$ a.e.
\end{itemize}
\end{proposition}
}

The variational inequality \eqref{stop} can be equivalently written as the inclusion
$\dot z(t) + \partial I_Z(z(t)) \ni \dot u(t)$.
This enables us to rewrite the differential inclusions \eqref{phase}
and \eqref{fati} for the phase variable $\chi$ and fatigue variable $m$
with a choice $\chi^0(x) \in [0,1]$, $m^0(x) \ge 0$ of initial conditions in the form
\begin{eqnarray}
\label{ph3.4}
\chi(x,t) &=& \stop_{[0,1]} [\chi^0(x), A(x,\cdot)](t),\\ \label{ph3.5}
m(x,t) &=& \stop_{[0,\infty)} [m^0(x), S(x,\cdot)](t),
\end{eqnarray}
where
\begin{eqnarray}
\label{ph3.2}
A(x,t) &:=& \int_0^t\frac{1}{\gamma}\left(\frac{L}{\thre}(\theta-\thre) \right)(x,\tau)\,\dd\tau,\\
\label{ph3.6}
S(x,t) &:=& \int_0^t\left( -h(\chi_t(\tau))
+ \int_0^1 \la(x-y) D[m,\ve](y,\tau)\dd y \right)(x,\tau)\,\dd\tau.
\end{eqnarray}
The advantage of this representation is that now, $\chi$ and $m$ are defined by equations
involving, by virtue of Proposition \ref{3p2}, only operators that are Lipschitz continuous
in $C[0,T]$ and in $W^{1,1}(0,T)$.

The variational inequality \eqref{stop} is also used to model single-yield elastoplasticity.
In this case, the constraint $Z=[-r,r]$ is the admissible stress domain, the input $u = \ve$
is the strain, and the output $z = \sigma_r:= \stop_r[\ve]$ is the stress.
We can rewrite \eqref{stop} equivalently in ``energetic'' form
\be{pi2a}
\dot \ve(t) \sigma_r(t) = \frac{\dd}{\dd t} \left( \frac12 \sigma_r^2(t)\right) + r|\dot\xi(t)|.
\ee
Indeed, $\dot\ve(t) \sigma_r(t)$ is the power supplied to the system, part of it is
used for the increase of the potential $\frac12 \sigma_r^2(t)$, and the rest $r|\dot\xi(t)|$
is dissipated.

The Prandtl-Ishlinskii model is constructed as a linear combination
of stops with all possible yield points
$r>0$. Here, given a measurable function $\varphi: [0, \infty) \times (0,\infty) \to [0,\infty)$
satisfying Hypothesis \ref{h1}\,(i) below,
we define the fatigue dependent
Prandtl-Ishlinskii operator $P: (W^{1,1}(0,T))^2 \to W^{1,1}(0,T)$ by the integral
\be{pi1}
P[m, \ve](t) = \int_0^\infty \varphi(m(t), r)\; \stop_r[\ve](t)\dd r\,.
\ee
Eq.~\eqref{pi2a} enables us to establish the energy balance for the Prandtl-Ishlinskii operator
\eqref{pi1}. Indeed, if we define the Prandtl-Ishlinskii potential
\be{pot}
V[m, \ve](t) = \frac12\int_0^\infty \varphi(m, r) \stop_r^2[\ve](t)\dd r\,,
\ee
and the dissipation operator
\be{diss}
D[m, \ve](t) = \int_0^\infty r \varphi(m, r) |\play_r[\ve]_t(t)|\dd r\,,
\ee
we can write the Prandtl-Ishlinskii energy balance in the form
\be{pien}
\dot\ve(t) P[m, \ve](t) = \frac{\dd}{\dd t} V[m, \ve](t) + D[m, \ve](t)
- \frac12 m_t \int_0^\infty \varphi_m(m, r) \stop_r^2[\ve](t)\dd r\quad \ale
\ee
As a consequence of Proposition \ref{3p2}\,(iv), we have
\begin{equation}
\label{D-ep}
D[m, \ve](t) \le |\dot\ve(t)| \int_0^\infty r\varphi(m, r)\dd r\,.
\end{equation}


\section{Statement of the problem}\label{stat}

For any $T>0$, we denote $\Omega_T := (0,1)\times (0,T)$, $u(x,t) = \int_0^t \sigma(x,\tau)\dd\tau$,
$f(x,t) = \int_0^t F(x,\tau)\dd\tau + \rho w_t(x,0) - \alpha w_{xxt}(x,0)$.
We rewrite the equations \eqref{b3}, \eqref{motion}, \eqref{b7}, \eqref{phase}, \eqref{fati}
as the system
\bearr
u_t &=& B w_{xx}+ P[m,w_{xx}] + \nu w_{xxt} -\beta(\theta-\thre),\label{s1}\\
\rho w_t- \alpha w_{xxt}&=&-u_{xx} + f(x,t),\label{s2}\\
c\theta_t- \kappa\theta_{xx}&=& D[m, w_{xx}] + \nu w_{xxt}^2- \beta\theta w_{xxt}
-\frac{1}{2}m_t\int^{\infty}_0\!\!{\varphi_m(m,r)\stop_r^2[w_{xx}]\dd r} \nonumber\\
&& -L\chi_t + g(\theta,x,t), \label{s3}\\
-\gamma\chi_t&\in&\partial I_{[0,1]}(\chi)-\frac{L}{\thre}(\theta-\thre),\label{s4}\\
m_t &\in & {- \partial I_{[0,\infty)}(m)} -h(\chi_t)+
\ipi \la(x-y) D[m, w_{xx}](y,t) \dd y,\label{s5}
\eearr
for unknown functions $u, w, \theta, m, \chi$  with initial and boundary conditions
\be{ini}
\left.
\barr{rl}
w(x,0) &= u(x,0) = 0\,,\\
{m(x,0)} &{= m^0(x) = 0}\,,\\
\theta(x,0) &= \theta^0(x),\,\\
\chi(x,0) &=\chi^0(x),
\earr
\right\}
\ee
\be{bou}
\left.
\barr{rl}
w(0,t) &= u(0,t) = w(1,t) = u(1,t) = 0\,,\\
\theta_x(0,t) &= \theta_x(1,t) = 0\,.\\
\earr
\right\}
\ee
The zero initial conditions for $w$ and $m$ are motivated by the fact that
it is difficult to determine the initial degree of fatigue for a material with unknown
loading history, and the most transparent hypothesis consists in assuming that
no deformation (and therefore no fatigue) has taken place prior to the time $t=0$.

The data are required to fulfill the following hypotheses:

\begin{hypothesis}\label{h1} \hfill
\begin{itemize}
\item[{\rm (i)}] $P$ is a Prandtl-Ishlinskii operator {\eqref{pi1} with distribution
function $\varphi: [0,\infty)\times (0,\infty) \to [0,\infty)$, locally Lipschitz
continuous in the first variable, and there exist
$\tilde{\varphi},\varphi^* \in L^1(0, \infty)$ such that
$\varphi(m,r) \le \tilde{\varphi}(r)$, $0 \le - \varphi_m(m,r) \le \varphi^*(r)$,
{$|\varphi_{mm}(m,r)| \le \varphi^*(r)$ a.e., with} 
${\tilde{M} :=} \int_0^{\infty} r\tilde{ \varphi }(r) \dd r < \infty$,
$M:= \int_0^{\infty} r^2 { \varphi^* }(r) \dd r < \infty$.}
\item[{\rm (ii)}] {$B, \nu, \beta, \thre, \rho, \alpha, c, \kappa, L, \gamma$ are given positive constants.}
\item[{\rm (iii)}] $\la: \real \to [0,\infty)$ is a $C^1$ function with compact support,
${\La}:= \max\{\la(x) + |\la'(x)|\,, \ x \in \real\}$.
\item[{\rm (iv)}] $f\in L^2(\Omega_T)$ is a given function for some fixed $T>0$, such that
$f_{tt}, f_{xt} \in L^2(\Omega_T)$.
\item[{\rm (v)}] $\theta^0 \in L^\infty(0,1)$ and $\chi^0\in W^{1,2}(\Omega)$ are
such that $\theta^0 \ge \theta_* > 0$, $\theta^0_{x} \in L^2(0,1)$, $\chi^0(x) \in [0,1]$
for all $x \in [0,1]$.
\item[{\rm (vi)}]   $h: \real \to [0,\infty)$ is a nondecreasing  Lipschitz continuous function
such that $h(z) \leq bz^2,$ $0 \le h'(z) \le a \ale$ for $z \in \real$, and $a,b$
are positive constants such that $bM \le \gamma$, {where $M$ is as in  (i) and $\gamma$ is the relaxation coefficient from  \eqref{s4}}.
\item[{\rm (vii)}] {$g:[0,\infty)\times \Omega_T \to \real$}
is a Carath\'eodory function and 
 $g_0(x,t) := g(0,x,t) \ge 0$,  for $g_0 \in L^2(\Omega_T)$ and  $|g_\theta(\theta,x,t)| \le g_1$ a.e. with  $g_1$  a constant.
\end{itemize}
\end{hypothesis}

The assumption that $\varphi(m,r)$ decreases with
increasing fatigue $m$ {corresponds to the observation that the stiffness} of the material
decreases with increasing fatigue. Also the assumption that $g_0(x,t) \ge 0$ makes sense.
Note that $g$ is the heat source density, so that at zero temperature,
we cannot remove heat from the system. 

We now check that regular solutions of \eqref{s1}--\eqref{bou} satisfy \eqref{b7}--\eqref{b8}
with $\ve = w_{xx}$ and $\sigma$ given by \eqref{b3}, which implies the thermodynamic consistency of the system.
Indeed, we have by \eqref{s3} and \eqref{pien}
$$
\UU(\ve, \theta, \chi)_t +  q_x - \ve_t \sigma = g
$$
and by \eqref{b5}, \eqref{s3}, and \eqref{s4}
$$
\SS(\ve, \theta, \chi)_t +  \left(\frac{q}{\theta}\right)_x - \frac{g}{\theta} = 
\frac{\kappa\theta_x^2}{\theta^2} + \frac{\nu\ve_t^2}{\theta}
+\frac{1}{\theta} \left(D[m, \ve](t)
- \frac12 m_t \int_0^\infty \varphi_m(m, r) \stop_r^2[\ve](t)\dd r + \gamma\chi_t^2\right).
$$
By Hypothesis \ref{h1} (i) and \eqref{s5},
$$
D[m, \ve](t)
- \frac12 m_t \int_0^\infty \varphi_m(m, r) \stop_r^2[\ve](t)\dd r + \gamma\chi_t^2
\ge \gamma\chi_t^2 - \frac{M}{2} h(\chi_t),
$$
hence, by Hypothesis \ref{h1}\,(vi)
$$
\SS(\ve, \theta, \chi)_t +  \left(\frac{q}{\theta}\right)_x \ge 
\frac{\kappa\theta_x^2}{\theta^2} + \frac{\nu\ve_t^2}{\theta} + \frac{\gamma\chi_t^2}{2\theta}
\ge 0,
$$
provided we check that the absolute temperature $\theta$ stays positive.
In Subsection \ref{loca}, we will find a positive lower bound for the discrete approximations
of the temperature, which is independent of the discretization parameter, and therefore is preserved
in the limit and implies the positivity of the temperature.

The main result of this paper reads as follows.

\begin{theorem}\label{t2}
Let Hypothesis \ref{h1} hold. Then there exists   a 
unique solution to the  system \eqref{s1}--\eqref{bou} in $\Omega_T$ {such that
$\theta(x,t) > 0$ for all $(x,t) \in \Omega_T$, and} with the regularity
\begin{itemize}
\item $w_{xxxt}, w_{xxtt}, \theta_t, \theta_{xx}, u_{tt}, u_{xxt} \in L^2(\Omega_T)$,
\item $m_t, {\chi_t} \in L^\infty(\Omega_T)$.
\end{itemize}
\end{theorem}


\section{Function spaces, interpolation}


Let $p,q,s \in [1, \infty]$ be such that $q > s$, and let $|\cdot|_p$ denote the norm in $L^p(0,1)$,  $||\cdot||_p$ the norm in $L^p(\Omega_T).$
 
 The
{\it Gagliardo-Nirenberg inequality\/} states that there exists a constant $C>0$
such that for every $v \in W^{1,p}(0,1)$ we have
\be{gag}
|v|_q \le C\left(|v|_s + |v|_s^{1-\vr} |v'|_p^\vr\right)), \quad \vr = \frac{\frac {1}{s} - \frac{1}{q}}
{1 + \frac {1}{s} - \frac{1}{p}}\,.
\ee
In fact, \eqref{gag} is straightforward: If we introduce an auxiliary parameter
$r = 1+ s(1-\frac{1}{p})$ and use the chain rule
$\frac{\dd}{\dd x} |v(x)|^r \le r |v(x)|^{r-1} |v'(x)|$ a.e., we obtain from H\"older's inequality
the estimate
$$
|v|_\infty \le |v|_r + C |v|_s^{1-(1/r)} |v'|_p^{1/r}.
$$
Combined with the obvious interpolation inequality $|v|_h \le |v|_\infty^{1-(s/h)}|v|_s^{s/h}$
for $h=q$ (and for $h=r$, if $r>s$), this yields \eqref{gag}.

Let now $\bfv = (v_0, v_1, \dots, v_n)$ be a vector, and let us denote
\be{vect}
|\bfv|_p = \left(\frac{1}{n} \sum_{k=0}^n |v_k|^p\right)^{1/p}\,, \quad
|\bfD\bfv|_p = \left(n^{p-1}\sumk |v_k- v_{k-1}|^p\right)^{1/p}.
\ee
The discrete counterpart of \eqref{gag} reads
\be{gagd}
|\bfv|_q \le C\left(|\bfv|_s + |\bfv|_s^{1-\vr} |\bfD\bfv|_p^\vr\right),
\quad \vr = \frac{\frac {1}{s} - \frac{1}{q}}{1 + \frac {1}{s} - \frac{1}{p}}\,,
\ee
and can be easily derived from \eqref{gag} by defining $v$ e.g. as equidistant
piecewise linear interpolations of $v_k$.


\section{Formal estimates}\label{form}

In order to explain the estimation technique, we first proceed formally, assuming that the positivity
of temperature is already established.
For the sake of simplicity we set from now on
\begin{eqnarray}\label{a}
&& \mathcal{K}[m, w_{xx}](x,t) := - \frac{1}{2} \int_0^{\infty} \varphi_m(m, r) \stop^2_r[w_{xx}] (x,t)\, \dd r\\\label{b}
&& \mathcal{D}[m, w_{xx}](x,t) := \ipi \la(x-y) D[m, w_{xx}](y,t) \dd y.
\end{eqnarray}
 Due to the fact that, by Definition \ref{dpi1} we have $|s_r[w_{xx}](t)| \le r$ and from Hypothesis \ref{h1}\,(i) we deduce
\begin{equation}
\label{mathcalK}
0 \leq \mathcal{K}[m, w_{xx}] \le \frac{1}{2} \int_0^{\infty} r^2 \varphi^*(r) \dd r = \frac{M}{2}.
\end{equation}
Finally, due to Hypothesis \ref{h1}  (i) {and (iii) and \eqref{D-ep}}, we have
\begin{equation}
\label{mathcalD}
0 \leq \mathcal{D}[m, w_{xx}](x,t) \le \, \La \tilde{M} \, |w_{xxt}(t)|_1.
\end{equation}
We will denote  in the sequel by  $C$ any constant possibly depending on
{the constants in Hypothesis \ref{h1}} and on the initial data of the problem.

\subsection{The energy estimate}

We {multiply} \eqref{s1} by $w_{xxt}$, obtaining
\begin{equation}
\label{test1}
- u_t \, w_{xxt} + B w_{xx} \, w_{xxt} + P[m, w_{xx}] w_{xxt} + \nu w_{xxt}^2
- \beta (\theta - \thre) w_{xxt} = 0;
\end{equation}
we differentiate \eqref{s2} in time and multiply by $w_t$, getting
\begin{equation}
\label{test2}
\rho w_{tt} w_t - \alpha w_{xxtt} w_t + u_{xxt} w_t - f_t w_t = 0,
\end{equation}
and finally we  sum up \eqref{test1}, \eqref{test2} and \eqref{s3}, {all } integrated in space. 
The first term in \eqref{test1} simplifies with the third term in \eqref{test2} due to integration by parts and our choice of boundary conditions; moreover the viscosity terms cancel out and also a term  $\beta \,\theta$. Concerning the term with hysteresis, using \eqref{pien} we deduce 
\begin{equation}
\label{EB-hyst}
P[m, w_{xx}] w_{xxt} = \frac{\dd}{\dd t} V[m, w_{xx}] + D[m, w_{xx}] + m_t \, \mathcal{K}[m, w_{xx}]
\end{equation}
and thus in the sum of \eqref{test1}, \eqref{test2} and \eqref{s3} what remains is just the term containing $V$. More precisely we have the  {\sl energy balance}
\begin{equation}
\label{balance}
\frac{\dd}{\dd t} \int_0^1 \left (\frac{1}{2} B w_{xx}^2 + V[m, w_{xx}] + \beta \thre w_{xx} + \frac{1}{2} \rho w_t^2 + \frac{1}{2}\alpha  w_{xt}^2 + c \theta + L \chi \right ) \dd x = \int_0^1 (f_t \, w_t \,+ g )
\dd x,
\end{equation}
and Gronwall's argument together with Hypothesis \ref{h1} (iv) and (vii) gives the estimate
\be{en-reg}
\forall t\in [0,T]:\  |w_{xx}(t)|_2 + |w_t(t)|_2 + |w_{xt}(t)|_2 + |\theta(t)|_1  \le C.
\ee


\subsection{The Dafermos estimate}

\label{dafer}

We test the equation for the temperature \eqref{s3} by $\theta^{-1/3}$ and obtain, using { notations} \eqref{a} and \eqref{b}, that
\begin{eqnarray}
0 &=& \int_0^1  - c \theta_t \theta^{-1/3} \dd x + \int_0^1 \kappa \theta_{xx} \theta^{- 1/3}  \dd x + \int_0^1 m_t \mathcal{K}[m, w_{xx}] \theta^{- 1/3}  \dd x + \int_0^1 D[m, w_{xx}] \theta^{-1/3}  \dd x \nonumber \\ && + \int_0^1 \nu w_{xxt}^2 \theta^{-1/3} \dd x 
 - \int_0^1 \beta w_{xxt} \theta^{2/3}  \dd x - \int_0^1 L \chi_t \theta^{-1/3} \dd x + \int_0^1 g \, \theta^{-1/3} \dd x \nonumber\\ 
 &=:& T_1 + T_2 + T_3 + T_4 + T_5 + T_6 + T_7 + T_8. \label{daf-3}
\end{eqnarray}
{We keep the terms $T_5 =\nu  \int_0^1  w_{xxt}^2 \theta^{-1/3} \dd x$,
$T_6 = - \beta \int_0^1 w_{xxt} \theta^{2/3}\dd x$, and
\be{T2}
T_2 = \int_0^1 \kappa \theta_{xx} \theta^{- 1/3}  \dd x
= \frac{\kappa}{3} \int_0^1 \left [\theta^{- 2/3} \theta_x \right ]^2 \dd x =
3 \kappa \int_0^1  \left [{(\theta^{1/3})_x} \right ]^2 \dd x,
\ee
where we have integrated by parts and used the boundary conditions \eqref{bou}.
All the other terms will be estimated from below. First,
\[
T_1 = - c \int_0^1  \theta_t \theta^{-1/3} \dd x =
- \frac{3c}{2}\frac{\dd}{\dd t}\int_0^1 (\theta^{2/3}) \dd x.
\]
In the identity
\be{T3}
T_3 = \int_0^1 m_t\, \mathcal{K}[m, w_{xx}] \theta^{- 1/3} \dd x
\ee
we have $T_3 \ge 0$ whenever $m_t \ge 0$. On the other hand, if $m_t < 0$, then by \eqref{s5}, \eqref{b}, and Hypothesis \ref{h1}\,(vi) we have
\be{fa1}
m_t = -h(\chi_t) + \mathcal{D}[m, w_{xx}] \ge -a \chi_t.
\ee
Now the assumption $m_t < 0$ implies that $\chi_t >0$. Then, however,
by \eqref{s4}, we have that
\begin{equation}
\label{vachi}
\gamma \chi_t = \frac{L}{\thre} (\theta - \thre) \le \frac{L}{\thre} \theta \qquad {a.e.}
\end{equation}
Combining the above inequalities with \eqref{mathcalK}, we  obtain for $m_t < 0$ that
\be{T3a}
T_3 \ge  - \frac{M L a}{2\gamma \thre} \int_0^1 \theta^{2/3} \dd x.
\ee
We obviously have
\[
T_4 = \int_0^1 D[m, w_{xx}] \theta^{-1/3}  \dd x\ge 0.
\]
}
The term
\[
T_7 := - L \int_0^1 \chi_t \theta^{-1/3} \dd x
\]
can be treated in a similar way as the term $T_3$ and using \eqref{vachi} for $\chi_t \ne 0$
we get
\[
T_7 \ge - \frac{L^2}{\gamma \thre} \int_0^1\theta^{2/3}	\dd x.
\] 
Finally, we find a lower bound for $T_8$ by Hypothesis \ref{h1}\,(vii) as follows:
{\[ 
T_8 = {\int_0^1} g(\theta,x,t) \, \theta^{-1/3} \dd x  \ge 
{\int_0^1} (g(\theta,x,t) - g(0,x,t)) \, \theta^{-1/3} \dd x
\ge  - g_1 \, \int_0^1  \theta^{2/3} \dd x.
\]
}
Coming back to \eqref{daf-3}, integrating it {in time, we deduce}
\begin{eqnarray}
&&  3\kappa \int_0^t \int_0^1 \left [{(\theta^{1/3})_x} \right ]^2 \dd x \dd\tau
+ \nu \int_0^t \int_0^1 w_{xxt}^2 \theta^{-1/3} \dd x \dd\tau\le \frac{3c}{2} \int_0^1 \theta^{2/3}\dd x  \nonumber\\
&& + \,C_1 \int_0^t \int_0^1  \theta^{2/3}\dd x\, \dd \tau  + \beta \int_0^t \int_0^1 |w_{xxt}| \theta^{2/3} \dd x \,\dd \tau, \label{dafermos}
\end{eqnarray}
where we put
\[
C_1 := \frac{L}{\gamma \thre}  \left (\frac{a M}{2} + L \right ) + g_1.
\]

The first two terms on the right hand side of \eqref{dafermos} are bounded due to \eqref{en-reg}.  The last term we estimate by H\"older inequality as follows
{
\begin{eqnarray*}
&&\hspace{-16mm}\beta \int_0^t \int_0^1 |w_{xxt}| \theta^{2/3}\dd x\,\dd \tau
= \beta \int_0^t \int_0^1 \theta^{5/6} \theta^{-1/6} |w_{xxt}| \dd x\,\dd \tau \\
&\le& \beta \int_0^t \left [ \left ( \int_0^1 \theta^{5/3}\dd x \right )^{1/2}
\left (\int_0^1 w_{xxt}^2  \theta^{- 1/3}\dd x \right )^{1/2}  \right  ]\dd \tau \\
&\le& \frac{\nu}{2} \int_0^t \int_0^1 w_{xxt}^2 \theta^{-1/3} \dd x\,\dd \tau + \frac{\beta^2}{2 \nu} \int_0^t \int_0^1 \theta^{5/3}\dd x\,\dd \tau, 
\end{eqnarray*}
and \eqref{dafermos} yields that}
\begin{equation}
\label{daf-1}  \int_0^t \int_0^1 \left [ {(\theta^{1/3})_x} \right ]^2 \dd x\,\dd \tau +  \int_0^t \int_0^1 w_{xxt}^2 \theta^{-1/3} \dd x\,\dd \tau \le \, C \left (1 +  \int_0^t \int_0^1 \theta^{5/3}\dd x\,\dd \tau \right). 
\end{equation}
Now we apply \eqref{gag} with {$v = \theta^{1/3}$, $s = 3$, $q = 5$, $p = 2$, $\rho = {4}/{25}$
and notice that}
\[
|\theta^{1/3}|_3 = \left (\int_0^1 (\theta^{1/3})^3\dd x \right )^{1/3}  \le C
\]
due to \eqref{en-reg}. We therefore have
\[
\left (\int_0^1 (\theta^{1/3})^5 \dd x\right )^{1/5} = |\theta^{1/3}|_5 \le \, C \left (|\theta^{1/3}|_3 + |\theta^{1/3}|_3^{21/25} \, \left|{(\theta^{1/3})_x}\right|_2^{4/25} \right )
\]

so that
\[
\int_0^t \int_0^1 \theta^{5/3} \dd x \dd \tau
\le C \left (1 + \int_0^t \left [\int_0^1 \left [{(\theta^{1/3})_x}\right ]^2 \dd x\right]^{2/5}  \dd\tau \right ).
\]
Combining this last estimate with \eqref{daf-1}, we deduce 
\begin{equation}\label{888}
\int_0^t \int_0^1 \left [{(\theta^{1/3})_x} \right ]^2 \dd x\dd\tau
+  \int_0^t \int_0^1 w_{xxt}^2 \theta^{-1/3}\dd x \dd \tau \le  C.
\end{equation}
Applying {again} the Gagliardo-Nirenberg inequality with the choices
$v = \theta^{1/3}$, $s = 3$, $q = 8$, $p = 2$, $\rho = {1}/{4},$ we obtain that
\[
\int_0^1 \theta^{8/3} \dd x  = |\theta^{1/3}|_8^8
\le C \left (|\theta^{1/3}|_3 + |\theta^{1/3}|_3^{3/4}
\left|{(\theta^{1/3})_x}\right|_2^{1/4}\right)^8
\]
and this   {after space integration,} together with \eqref{en-reg} and \eqref{888}  brings the estimate
\begin{equation}
\label{reg-2-theta}
\|\theta\|_{8/3} = \int_0^T \int_0^1 \theta^{8/3}\dd x\dd t \le C.
\end{equation}
}
To derive a further estimate, we sum again \eqref{test1} and \eqref{test2}, and obtain
\begin{eqnarray*}
&&\hspace{-10mm} {\frac{\dd}{\dd t} \int_0^1 \left (\frac{B}{2} w_{xx}^2 + \beta \thre w_{xx}
+ {\frac{\rho}{2} w_t^2} + \frac{\alpha}{2} w_{xt}^2 \right) \dd x
+ \int_0^1 \nu w_{xxt}^2  \dd x} \\
&=& \int_0^1 {\left(\beta\theta w_{xxt} - P[m, w_{xx}]\, w_{xxt} + f_t w_t\right)} \dd x.
\end{eqnarray*}
{We estimate the first term on the right hand side using the inequality
$\beta\theta w_{xxt} \le \frac{\beta^2}{2\nu} \theta^2 + \frac{\nu}{2} w_{xxt}^2$
and the previous estimate \eqref{reg-2-theta}. {In the second term, $P[m, w_{xx}]$
is bounded by Hypothesis \ref{h1}\,(i), and the} third term is handled using Hypothesis \ref{h1}\,(iv).
This finally gives the additional estimate 
\begin{equation}
\label{epst2}
\|w_{xxt}\|_2 \le C.
\end{equation}
}


\subsection{Higher order estimates}

We differentiate \eqref{s1}  in space, obtaining 
\begin{equation}
\label{s1space}
u_{xt} = B w_{xxx} + P[m, w_{xx}]_x + \nu w_{xxxt} - \beta \theta_x.
\end{equation}
We integrate by parts in space, recalling that if $u(0,t) = u(1,t) = 0$ {by \eqref{bou},
then $u_t(0,t) = u_t(1,t) = 0$. We deduce by \eqref{s1} and \eqref{s2}}
\begin{eqnarray}\nonumber
&& \int_0^1 (B w_{xxx} + P[m, w_{xx}]_x + \nu w_{xxxt} - \beta \theta_x)^2 \dd x \stackrel{\eqref{s1space}}{=} \int_0^1 u_{xt}^2  \dd x = \int_0^1 u_t (- u_{xxt}) \dd x\\ \label{integ}
&& {=} \int_0^1 (\rho w_{tt} - \alpha w_{xxtt} - f_t) \,(B w_{xx} + P[m, w_{xx}] + \nu w_{xxt} - \beta (\theta - \thre)) \,  \dd x.
\end{eqnarray}
This brings
\begin{eqnarray} \nonumber
&& \int_0^1 (\nu w_{xxxt} + B w_{xxx})^2 \dd x +  \nu \rho \int_0^1 w_{xt} w_{xtt} \dd x
+ \nu \alpha \int_0^1  w_{xxt} w_{xxtt} \dd x \\ \nonumber
&\le & {C \int_0^1 (P[m, w_{xx}]_x^2 + \theta_x^2) \dd x
- \int_0^1 f_t (B w_{xx} + P[m, w_{xx}] + \nu w_{xxt} - \beta (\theta - \thre)) \dd x}\\ \nonumber
&& +
\frac{\dd}{\dd t} \int_0^1(\rho w_{t}-\alpha w_{xxt})\, (B w_{xx} + P[m, w_{xx}] -\beta(\theta - \thre))\dd x 
\\ \label{hor1}
&& - \int_0^1 (\rho w_{t} - \alpha w_{xxt}) \, (B w_{xxt} + P[m, w_{xx}]_t  - \beta \theta_t) \dd x.
\end{eqnarray}
First of all, using Hypothesis \ref{h1} (i) and (iv), \eqref{en-reg}, \eqref{reg-2-theta} and \eqref{epst2} we estimate
\be{rhs1}
{\int_0^T \int_0^1 |f_t| |B w_{xx}+P[m, w_{xx}]+\nu w_{xxt}-\beta(\theta-\thre)|\dd x\dd t} \le C.
\ee
{Furthermore, by \eqref{en-reg} and Hypothesis \ref{h1} (i) 
we have
\be{rhs2}
\int_0^1(\rho w_{t}-\alpha w_{xxt})\, (B w_{xx} + P[m, w_{xx}] -\beta(\theta - \thre))\dd x
\le C(1 + |w_{xxt}|_2)(1 + |\theta|_2).
\ee
Note that by \eqref{reg-2-theta} we have
\begin{equation}\label{th}
{|\theta|_2^2 - |\theta^0|_2^2 = 2\int_0^t\int_0^1 \theta \theta_t \dd x\dd\tau \le
C\left(\int_0^t |\theta_t|_2^2 \dd \tau\right)^{1/2},}
\end{equation}
and \eqref{rhs2} yields that
\bearr \nonumber
&&\hspace{-14mm}\int_0^1(\rho w_{t}-\alpha w_{xxt})\, (B w_{xx} + P[m, w_{xx}] -\beta(\theta - \thre))
\dd x
\\ \nonumber
&\le& C(1 + |w_{xxt}|_2)\left(1 + \int_0^t |\theta_t|_2^2 \dd \tau\right)^{1/4}\\ \label{rhs3}
&\le& \frac{\nu\alpha}{4}|w_{xxt}|_2^2 + C\left(1 + \int_0^t |\theta_t|_2^2 \dd \tau\right)^{1/2}.
\eearr
Finally, still by \eqref{en-reg} and \eqref{epst2},
\bearr \nonumber
&&\hspace{-14mm}\left|\int_0^t\int_0^1 (\rho w_{t} - \alpha w_{xxt}) \, (B w_{xxt} + P[m, w_{xx}]_t
- \beta \theta_t)  \dd x \dd\tau\right|\\ \label{rhs4}
 &\le& C\left(1 + \int_0^t |P[m, w_{xx}]_t|_2^2
+|\theta_t|_2^2 \dd \tau\right)^{1/2}.
\eearr
Note that we have for a.e. $(x,t) \in \Omega_T$ that
$$
{\left| P[m, w_{xx}]_t(x,t) \right| \le \left| m_t \int_0^{\infty} \varphi_m(m,r) \stop_r[w_{xx}] \,\dd r\right| + \left| \int_0^{\infty} \varphi(m,r) (\stop_r[w_{xx}])_t \dd r \right|.}
$$
By \eqref{stop} we have $|\stop_r[w_{xx}]| \le r$, and Proposition \ref{3p2}\,(iv) yields
\be{rhs0}
|\stop_r[w_{xx}]_t| \le |w_{xxt}|,\
|m_t| \le |- h(\chi_t) + \mathcal{D}[m, w_{xx}]| \le C (|\chi_t| + |w_{xxt}|_1), \ 
|\chi_t| \le C (1 + \theta)\ \ale
\ee
{}From Hypothesis \ref{h1}\,(i) we thus obtain the pointwise bound
\be{pt1}
{\left| P[m, w_{xx}]_t \right| \le  C (1 + \theta + |w_{xxt}|_1) \ \ale,}
\ee
and from \eqref{rhs4}, using {\eqref{reg-2-theta}} and \eqref{epst2}  we conclude that
\be{rhs5}
\left|\int_0^t\int_0^1 (\rho w_{t} - \alpha w_{xxt}) \, (B w_{xxt} + P[m, w_{xx}]_t
- \beta \theta_t)  \dd x \dd\tau\right|
\le C\left(1 + \int_0^t|\theta_t|_2^2 \dd \tau\right)^{1/2}.
\ee
We now integrate \eqref{hor1} from $0$ to $t$, {for $t \in (0,T)$}. 
We combine \eqref{rhs1}, \eqref{rhs3}, \eqref{rhs5}, and \eqref{eini1} with \eqref{hor1}
integrated in time
and obtain
\begin{eqnarray} \nonumber
&&\hspace{-14mm} \int_0^1 \left (w_{xt}^2 + w_{xxt}^2 + w_{xxx}^2\right ) \dd x + \int_0^t \int_0^1 (w^2_{xxx} + w^2_{xxxt})  \dd x  \dd \tau
\\ \label{hor2}
&\le & C \left (1 + \int_0^t \int_0^1 (P[m, w_{xx}]^2_x + \theta_x^2)\dd x \dd \tau 
+ \left(\int_0^t |\theta_t|_2^2 \dd \tau\right)^{1/2}\right).
\end{eqnarray}
Here we  had to   estimate the initial values
$$
\int_0^1 (w_{xt}^2 + w_{xxt}^2 + w_{xxx}^2)(x,0)\dd x,
$$ which can be done as follows:
We have by \eqref{ini} $w_{xxx}(x,0) = 0$. Eq.~\eqref{s2} for $t=0$ reads
$$
\rho w_{t}(x,0) - \alpha w_{xxt}(x,0) = f(x,0),
$$
and testing this identity by $w_{xxt}(x,0)$ we see that
\be{eini1}
\int_0^1 (w_{xt}^2 + w_{xxt}^2 + w_{xxx}^2)(x,0)\dd x \le C.
\ee

Finally,} we deal with the term
\be{rhs6}
{P[m, w_{xx}]_x(x,t)  = \int_0^{\infty} \left(\varphi_m(m,r) \;m_x\; \stop_r[w_{xx}]\right)(x,t) \dd r + \int_0^{\infty} \left(\varphi(m, r) \stop_r[w_{xx}]_x\right)(x,t) \dd r.}
\ee
For all $x, h$, and $t$, we have by Proposition \ref{3p2} (iii)
\[
|\stop_r[w_{xx}](x+h, t) - \stop_r[w_{xx}](x,t)| \le  2 \max_{\tau \in [0,t]} |{w_{xx}(x+h, \tau)} - w_{xx}(x, \tau)|,
\]
which implies
\[
|\stop_r[w_{xx}]_x(x,t)| \le \, 2 \max_{\tau \in [0,t]} |w_{xxx}(x, \tau)| \ \ale
\]
By \eqref{ph3.4}, \eqref{ph3.2}  and Proposition \ref{3p2}\,(ii), we have
\[
\int_0^t |\chi_t(x+h, \tau) - \chi_t(x,\tau)| \dd \tau \le \, C \left (|\chi(x+h, 0) - \chi(x,0)| + \int_0^t |\theta(x+h, \tau)  - \theta(x,\tau)| \dd \tau \right ),
\]
which entails {for a.e. $x \in (0,1)$ that}
\[
\int_0^t|\chi_{xt}(x, \tau)| \dd \tau \le \, C \left (|\chi^0_x(x)| + \int_0^t |\theta_x(x, \tau)| \dd \tau \right ),
\]
and in {a similar} way we obtain from \eqref{ph3.5}, {\eqref{ph3.6}, Hypothesis \ref{h1} points (i) and (vi), \eqref{D-ep}}, Proposition \ref{3p2}\,(ii), and \eqref{epst2}
that
\bears
|m_x(x,t)| &\le& \int_0^t |m_{xt}(x,\tau)|\dd\tau \le C \left( \int_0^t|\chi_{xt}(x,\tau)|\dd\tau
+ \int_0^t |w_{xxt}|_1(\tau)\dd\tau\right)\\
&\le& C \left(1 + \int_0^t|\chi_{xt}(x,\tau)|\dd\tau\right),
\eears
{where we also used  that $m(x,0) = 0$.}
Therefore, from \eqref{rhs6}, using Hypothesis \ref{h1}\,{(i) and (iii)} and \eqref{ini}, we get
{for a.e. $(x,t) \in \Omega_T$ that}
\begin{eqnarray}\nonumber
&&\hspace{-14mm} |P[m, w_{xx}]_x(x,t)| \le \, C \left (1 + |\chi_x^0(x)| + \int_0^t |\theta_x(x, \tau)| \dd \tau + \max_{\tau \in [0,t]} |w_{xxx}(x,\tau)| \right )\\ \label{hor3}
&\le& C {\left (1 + |\chi_x^0(x)| + \int_0^t |\theta_x(x, \tau)| \dd \tau
+ \int_0^t |w_{xxxt}(x,\tau)| \dd \tau \right).}
\end{eqnarray}
{Combining} \eqref{hor2} and \eqref{hor3} and applying Gronwall's lemma we deduce
\begin{equation} \label{laste}
| w_{xxt}(t)|_2^2 +| w_{xxx}(t)|_2^2 + 
|| w^2_{xxxt} ||_2^2
\le  C \left(1 +  ||  \theta_x||_2^2  +   || \theta_t ||_2 \right),
\end{equation}
{where  we also used Hypothesis \ref{h1}  (v)  to estimate $|\chi^0_x(x)|$.}

It remains to estimate the $W^{1,2}$-norm of $\theta$ {both in space and time}. In the first step, we test \eqref{s3} by
$\theta$ and obtain, using \eqref{rhs0}, {\eqref{D-ep}, \eqref{epst2} and Hypothesis \ref{h1} points (i) and (vii)}, that
\bearr \nonumber
\frac{\dd}{\dd t} \int_0^1 \theta^2 \dd x + \int_0^1 \theta_x^2 \dd x
&\le& C {\int_0^1}\left ({|w_{xxt}|\theta + \theta |w_{xxt}|^2 + \theta^2 |w_{xxt}| + |m_t| \theta +  |\chi_t| \theta +   \theta |g |}\right ) \dd x\\ \label{rht1}
&\le& {C \left(1 + \int_0^1 \left( \theta^2
+ w_{xxt}^2 + \theta w_{xxt}^2 + \theta^2 |w_{xxt}|\right ) \dd x\right)}
\eearr
{By H\"older's inequality and \eqref{reg-2-theta}, we have
\[
\int_0^t\int_0^1 \theta^2 |w_{xxt}|\dd x\dd t \le
\|\theta\|_{8/3}^2 \left(\int_0^t\int_0^1w_{xxt}^4\dd x\dd t\right)^{1/4}\le
C\left(\int_0^t\int_0^1 w_{xxt}^4\dd x\dd t\right)^{1/4},
\]
and
\[
\int_0^t\int_0^1 \theta w_{xxt}^2\dd x\dd t \le
\|\theta\|_2^2 \left(\int_0^t\int_0^1 w_{xxt}^4\dd x\dd t\right)^{1/2}\le
C \left(\int_0^t\int_0^1 w_{xxt}^4\dd x\dd t\right)^{1/2}.
\]
Exploiting once more \eqref{reg-2-theta} {and \eqref{epst2}, we obtain finally integrating \eqref{rht1}
with respect to $t$ that}
\begin{equation}
\label{rht2}
 |\theta(t)|_2^2 + || \theta_x||_2^2   \le \, C \left (1 + || w_{xxt}||_4^2  \right ).
\end{equation}
On the other hand,  testing \eqref{s3} by $\theta_t$ we  deduce
\[
\int_0^1 \theta_t^2 \dd x + \frac{\dd}{\dd t} \int_0^1 \theta_x^2 \dd x \le C \, \int_0^1 (m_t^2 + w_{xxt}^2 + \chi_t^2) \dd x + \int_0^1 w_{xxt}^4 \dd x+ \int_0^1 \theta^2 w_{xxt}^2 \dd x
\]
and {by a similar argument as above} we obtain 
\begin{equation}
\label{rht3}
{\| \theta_t\|_2^2 + \sup_{t \in [0,T]} |\theta_x(t)|_2^2 \le C \left(1 + \|w_{xxt}\|_4^4 
+ \|\theta\|_4^4 \right).}
\end{equation}
We have by \eqref{reg-2-theta} that
\be{rht4}
\|\theta\|_4^4 \le \|\theta\|_{\infty}^{4/3} \|\theta\|_{8/3}^{8/3} \le C  \|\theta\|_{\infty}^{4/3}.
\ee
We now apply the Gagliardo-Nirenberg inequality \eqref{gag} with $q = \infty$, $s = 1$, $p = 2$,
and $\gamma = 2/3$, to deduce
\[
|\theta|_{\infty} \le C \left (|\theta|_1 + |\theta|_1^{1/3} |\theta_x|_2^{2/3}\right ).
\]
Using \eqref{en-reg} we obtain
\be{rht5}
\|\theta\|_{\infty} \le C \left (1 + \sup_{t \in [0,T]} |\theta_x|_2^{2/3} \right ).
\ee
{It follows from \eqref{rht3}, \eqref{rht4}, \eqref{rht5} that
\begin{equation}
\label{rht6}
\|\theta_t\|_2^2 + \sup_{t \in [0,T]} |\theta_x(t)|_2^2 \le C \left(1 + \|w_{xxt}\|_4^4 \right ), 
\end{equation}
which is what we were looking for. 
Coming back to  \eqref{laste}, using \eqref{rht2} and \eqref{rht6} we deduce in particular 
\be{last}
\sup_{t \in [0,T]} |w_{xxt}(t)|_2^2 + \|w_{xxxt}\|_2^2 \le C (1 + \|w_{xxt}\|_4^2).
\ee
We estimate the right hand side of \eqref{last} using the Gagliardo-Nirenberg inequality  \eqref{gag} with $q = 4$, $s = p = 2$, $\gamma = 1/4$,
and obtain
\[
|w_{xxt}(t)|_4 \le \, C \left (|w_{xxt}(t)|_2 + |w_{xxt}(t)|_2^{3/4} \, |w_{xxxt}(t)|_2^{1/4} \right),
\]
and this implies, by H\"older's inequality and \eqref{epst2} that
\begin{eqnarray} \nonumber
\|w_{xxt}\|_4^4 &\le& C \sup_{t \in [0,T]} |w_{xxt}(t)|_2^2 \left (||w_{xxt}||_2^2
+ \|w_{xxt}\|_2 \|w_{xxxt}\|_2\right)\\
&\le& C \sup_{t \in [0,T]} |w_{xxt}(t)|_2^2 \left (1 + \|w_{xxxt}\|_2\right ).\label{I}
\end{eqnarray}
}
Using this last estimate and coming back to \eqref{last} we get
\begin{eqnarray*}
||w_{xxt}||_4^2 &\le& C \sup_{t \in [0,T]} |w_{xxt}(t)|_2  \left (1 + ||w_{xxxt}||^{1/2}_2 \right ) \le \, \left (1 + ||w_{xxt}||_4^{3/2} \right ),
\end{eqnarray*}
{which enables us to conclude that
\be{rht7}
\|w_{xxt}\|_4 \le C,
\ee
and consequently by {\eqref{rht5}--\eqref{rht6}}
\begin{equation}
\label{rht8}
{\|\theta\|_\infty^2 + } \|\theta_t\|_2^2 + \sup_{t \in [0,T]} |\theta_x(t)|_2^2 \le C.
\end{equation}
Coming back to {\eqref{rhs0} and} \eqref{laste}, we deduce the following additional estimates
\be{rht9}
{\|m_t\|_\infty + \|\chi_t\|_\infty +} \sup_{t \in [0,T]}(|w_{xxt}(t)|_2
+ |w_{xxx}(t)|_2) + \|w_{xxxt}\|_2 \le C.
\ee
}
Finally, we differentiate \eqref{s1} in $t$ and test by
$w_{xxtt}$, differentiate \eqref{s2} twice in $t$ and test by $w_{tt}$, sum up the results,
eliminating the terms in $u_{tt}$ by integrating by parts. {We have to estimate
the initial values
$$
\int_0^1 (|w_{tt}^2(x,0)| + |w_{xtt}^2(x,0)|)\dd x.
$$
To do that, we proceed similarly as in \eqref{eini1}. We have by \eqref{ini} and \eqref{s1}
$$
u_t(x,0) = \nu w_{xxt}(x,0) - \beta(\theta^0(x) - \thre),
$$
hence
\be{eini2}
u_{xt}(x,0) = \nu w_{xxxt}(x,0) - \beta \theta^0_{x}(x).
\ee
On the other hand, by \eqref{s2},
\be{eini3}
\rho w_{xt}(x,0) - \alpha w_{xxxt}(x,0) = f_{x}(x,0),
\ee
and
\be{eini4}
\rho w_{tt}(x,0) - \alpha w_{xxtt}(x,0) = - u_{xxt}(x,0) + f_{t}(x,0).
\ee
We test \eqref{eini3} by $w_{xxxt}(x,0)$, use \eqref{eini1}, and obtain
\be{eini5}
\int_0^1 w_{xxxt}^2(x,0) \dd x \le C.
\ee
Hence, by \eqref{eini2} and Hypothesis \ref{h1} (v)
\be{eini6}
\int_0^1 u_{xt}^2(x,0) \dd x \le C.
\ee
Testing \eqref{eini4} by $w_{tt}(x,0)$ and integrating by parts we finally obtain 
\be{eini7}
\int_0^1 (\rho|w_{tt}^2(x,0)| + \alpha|w_{xtt}^2(x,0)|)\dd x \le \int_0^1 (|u_{xt}(x,0)|\, |w_{xtt}(x,0)|
+ |f_{t}(x,0)|\, |w_{tt}(x,0)|)\dd x,
\ee
which implies the desired estimate
\be{eini8}
\int_0^1 (|w_{tt}^2(x,0)| + |w_{xtt}^2(x,0)|)\dd x \le C.
\ee
This enables us to conclude that
\be{deri}
\sup_{t \in [0,T]}(|w_{tt}(t)|_2 + |w_{xtt}(t)|_2) + {\|w_{xxtt}\|_2} \le C.
\ee
}


\section{Approximation}\label{appr}

Here, we make rigorous the estimates derived formally in the previous section.
From now on, the values
of all physical constants are set to 1 for simplicity.

We choose an integer $n \in \nat$, and consider the space
discrete approximations of \eqref{s1}--\eqref{s4} for $k=1, \dots n-1$:
\begin{align}\label{ds1}
\dot u_k &= \ve_k + P[m_k, \ve_k] + \dot \ve_k - \theta_k + \thre\,,\\ \label{ds2}
\dot w_k - \dot\ve_k &= -n^2(u_{k+1} - 2 u_k + u_{k-1}) + f_k \,,\\ \label{ds3}
\ve_k &= n^2(w_{k+1} - 2 w_k + w_{k-1})\,,\\   \label{ds4}
\dot\theta_k &= n^2(\theta_{k+1} - 2 \theta_k + \theta_{k-1})
{ +\; \dot{m}_k\; \KK_k +D_k
+ \dot\ve_k^2 - \theta_k\dot\ve_k - \dot{\chi}_k + g_k(\theta_k,t)
\,,}\\ \label{ds5}
m_k  & = {\stop_{[0,\infty)}[0, S_k] ,\qquad S_k(t) = \int_0^t (-h(\dot{\chi}_k)+ \DD_k)(\tau)\dd \tau,}\\
\label{ds6}
{\chi}_k & = {\stop_{[0,1]}[\chi_k^0, A_k], \qquad A_k(t) = \int_0^t (\theta_k - \thre)(\tau) \dd \tau,} 
\end{align}
where 
\begin{eqnarray*}
&& {\KK_k(t) = - \frac{1}{2} \int_0^{\infty} \varphi_m(m_k(t), r) \stop_r^2[\varepsilon_k](t) \dd r
\in \left[0, \frac{M}{2}\right]\,,}\\
&& {D_k(t) = \int^{\infty}_0 \varphi(m_k(t),r)\, \stop_r[\varepsilon_k](t)
(\varepsilon_k -\stop_r[\varepsilon_k] )_t(t)\dd r \ge 0\,,}\\
&& {\DD_k(t) = \frac{1}{n} \sum_{j=1}^{n-1} \la_{j-k} D_j(t)\ge 0\,,}\\
&& {\la_i = \lambda(i/n)}\,,\\
&& f_k(t) = n \int_{(k-1)/n}^{k/n} f(x,t) \dd x \, ,\\[3.5mm]
&& g_k(\theta,t) = \left \{
\begin{array}{lll}
\!\!\!\!\!\! & \displaystyle n \int_{(k-1)/n}^{k/n} g(\theta,x,t) \dd x \qquad & \textnormal{for $\theta \ge 0$}\\[4.5mm]
\!\!\!\!\!\! & g_k(0,t) \qquad & \textnormal{for $\theta < 0.$}
\end{array}
\right.
\end{eqnarray*}
We prescribe initial conditions for $k=1, \dots, n-1$
\be{inid}
\left.
\barr{l}
w_k(0) = u_k(0) = 0 \,,\\[2mm]
\theta_k(0) = {\theta_k^0 :=} \theta^0(k/n)\,,\\[2mm]
m_k(0) = 0\,,\\[2mm]
\displaystyle \chi_k(0) = {\chi_k^0 :=} n \int_{(k-1)/n}^{k/n} \chi^0(x) \dd x\, ,
\earr
\right\}
\ee
and ``boundary conditions''
\be{bound}
\left.
\barr{l}
w_0 = w_n = u_0 = u_n = 0\,,\\
\theta_0 = \theta_1, \ \theta_n = \theta_{n-1}\,.
\earr
\right\}
\ee

This is a system of ODEs for $u_k, w_k, \theta_k $. We proceed as follows:

We claim   that \eqref{ds1}--\eqref{ds6} admits a $W^{1,\infty}$
solution in an interval $[0, T_n]$. First, denoting by $\bw$
the vector $(w_1, \dots, w_{n-1})$, and $\bve = (\ve_1, \dots, \ve_{n-1})$,
we have, by \eqref{ds3}, $-\bve = S \bw$ with a positive definite matrix $S$, which has the form 
\[
S = n^2 \left [
\begin{array}{cccccccccc}
2 & - 1 & 0 & 0 & 0 & \dots & 0 & 0 & 0 & 0\\
-1 & 2 & -1 & 0 & 0 & \dots & 0 & 0 & 0 & 0\\
0 & -1 & 2 & -1 & 0 & \dots & 0 & 0 & 0 & 0 \\
\vdots & \vdots & \vdots& \vdots& \vdots& \vdots& \vdots& \vdots &  \vdots& \vdots\\
0 & 0 & 0 & 0 & 0 & \dots & -1 & 2 & -1 & 0\\
0 & 0 & 0 & 0 & 0 & \dots & 0 & -1 & 2 & -1\\
0 & 0 & 0 & 0 & 0 & \dots & 0 &  0 & - 1 & 2
\end{array}
\right ]
\]
Hence, the left hand side of \eqref{ds2} reads $(I+S)\dot\bw$.
By \eqref{ds2}, $\dot\bve$  is itself a Lipschitz continuous mapping of 
$\mathbf{u} = (u_1, \dots, u_{n-1})$. Using Proposition \ref{3p2} (ii)  we see that \eqref{ds1}--\eqref{ds4}
can be considered as an ODE system in $u_k, w_k, \theta_k$, with a locally Lipschitz continuous and locally bounded
right hand side and the existence  and uniqueness of a local solution in an interval $[0,T_n]$
follows from the standard theory of ODEs, and the solution belongs to $W^{1,\infty}(0,T_n)$.

In the sequel, we will systematically use the ``summation by parts formula''
\be{parts}
\sumkj \xi_k(\eta_{k+1} - 2 \eta_k + \eta_{k-1}) + \sumk (\xi_k - \xi_{k-1})(\eta_k - \eta_{k-1})
= \xi_n(\eta_n - \eta_{n-1}) - \xi_0(\eta_1 - \eta_{0})
\ee
{for all vectors $(\xi_0, \dots, \xi_n)$, $(\eta_0, \dots, \eta_n)$.}


\subsection{Positivity of  the temperature}\label{loca}

In this subsection, we prove that $\theta_k$ remain positive in the whole
range of existence. As a first step, we test \eqref{ds4} by $-\theta_k^-$,
{where $\theta_k^-$ is the negative part of $\theta_k$.} 

We have
\[
- \frac{1}{n} \sum_{k=1}^{n-1} \dot{\theta}_k \theta_k^- = \frac{1}{2n} \frac{\dd}{\dd t} \sum_{k=1}^{n-1} (\theta_k^-)^2.
\]
On the other hand, by \eqref{parts} we also have
\[
- n \sum_{k=1}^{n-1} (\theta_{k+1} - 2 \theta_k + \theta_{k-1}) \theta_k^{-} = n \sum_{k=1}^n (\theta_k - \theta_{k-1}) (\theta_k^- - \theta_{k-1}^-) \le \, - n \sum_{n=1}^n (\theta_k^- - \theta_{k-1}^-)^2 \le 0,
\]
{Moreover, $D_k(t) \ge 0$ and $g_k(\theta, t) \ge 0$ for $\theta \le 0$ by Hypothesis \ref{h1}\,(vii),
hence
\[
- \left (D_k(t) + \dot{\ve}_k^2(t) + g_k(\theta_k,t)\right) \theta_k^- \le 0.
\]
}
Now we deal with the phase term. We have that
\begin{equation}
\label{phase-term}
\left.
\begin{array}{lll}
& \dot{\chi}_k(t) \theta_k^-(t) = 0 \qquad & \textnormal{if $\dot{\chi}_k(t) = 0$,}\\
& \dot{\chi}_k(t) \theta_k^-(t) = (\theta_k(t) - \thre) \theta_k^-(t) \le 0 \qquad & \textnormal{otherwise}.
\end{array}
\right \}
\end{equation}
Finally, {if $\dot{m}_k(t) \ne 0$, then
\[
-\dot{m}_k(t) \KK_k(t)  \theta_k^-(t) = \left (h(\dot{\chi}_k(t))- \DD_k(t)\right) \theta_k^-(t) \KK_k(t)
\le  h(\dot{\chi}_k(t)) \theta_k^-(t) \KK_k(t) \le 0
\]
by virtue of \eqref{phase-term}. Summarizing the above computations, we get}
\[
\frac{\dd}{\dd t} \frac{1}{2n} \sum_{k=1}^{n-1} (\theta_k^-)^2 \le \frac{1}{n} \sum_{k=1}^{n-1} (\theta_k^-)^2 \dot{\ve}_k \le \frac{K_{\ve, n}}{n} \sum_{k=1}^{n-1} (\theta_k^-)^2,
\]
where we put
\[
K_{\ve, n} := {\max \{|\dot{\ve}_k(t)|: k = 1, \dots, n-1, \,\,\, t \in [0, T_n]\}}
\]
and Gronwall's argument yields $\theta_k^-(t) = 0$ for all $k$ and $t \in [0,T_n]$.

{We now prove that in fact, $\theta_k(t)$ are} bounded away from $0$ for all $k$ and all $t \in [0,T_n]$.
First of all we notice that if $\dot{\chi}_k \neq 0$ then 
\[
- \dot{\chi}_k = - \theta_k + \thre \ge - \theta_k.
\]
On the other hand
\[
\dot{m}_k \ge - h(\dot{\chi}_k) \ge - a \dot{\chi}_k \ge - a \theta_k, 
\]
{hence
\[
\dot{m}_k \KK_k \ge - \frac{M a}{2} \theta_k.
\]
Using the estimates above together with 
\[
\dot{\ve}_k^2 - \theta_k \dot{\ve}_k \ge - \frac{1}{4} \theta_k^2
\]
{we obtain from \eqref{ds4} that
\[
\dot{\theta}_k - n^2 (\theta_{k+1} - 2 \theta_k + \theta_{k-1}) \ge - \psi(\theta_k),
\]
}
where we set
\[
\psi(z) := \frac{1}{4} z^2 + a \left (1 + \frac{M}{2} \right ) z.
\]
Let $p$ be the solution of the differential equation
$$
\dot p + \psi(p) = 0\,, \quad p(0) = {\theta_*},
$$
{with $\theta_*>0$ from Hypothesis \ref{h1}\,(v).} It is easy to check that 
\[
{p(t) = \frac{\mu \theta_* e^{- \mu t}}{\delta \theta_*(1-e^{- \mu t}) + \mu}\ ,  \qquad \text{with}
\qquad \delta = \frac{1}{4}, \qquad \mu = a \left (1 + \frac{M}{2} \right).} 
\]
Then
\be{pthetak}
(\dot p - \dot\theta_k) - n^2((p- \theta_{k+1}) - 2 (p-\theta_k) + (p-\theta_{k-1}))
+ {\psi (p) - \psi (\theta_k)}\le 0\,.
\ee
Testing \eqref{pthetak} by $(p-\theta_k)^+$ and using \eqref{parts}, we obtain
\be{e7a}
\frac12 \frac{\dd}{\dd t}\sumkj((p - \theta_k)^+)^2
+ n^2 \sumk ((p-\theta_k)^+ - (p-\theta_{k-1})^+)^2
+ {(\psi (p) - \psi (\theta_k))} (p - \theta_k)^+ \le 0\,,
\ee
hence, {as $\psi$ is nondecreasing for positive arguments,
\[
\sumkj((p - \theta_k)^+)^2(t) \le \sumkj((p - \theta_k)^+)^2(0) = 0, 
\]
so that $\theta_k(t) \ge p(t)  >0$
for all $k$ and all $t \in [0,T_n]$, which is the desired result.}


\subsection{Discrete energy estimate}\label{ener}

We test \eqref{ds1} by $\dot \ve_k$, differentiate \eqref{ds2} in $t$ and test by $\dot w_k$,
and sum up over $k=1, \dots n-1$. From \eqref{ds3}, with a repeated use of
\eqref{bound} and \eqref{parts}, we obtain
\be{es5x}
\frac{\dd}{\dd t} \left(\frac{1}{2n}\sumkj \dot w_k^2
+ \frac{n}{2} \sumk (\dot w_k - \dot w_{k-1})^2\right)
+ \frac{1}{n} \sumkj \dot \ve_k(\ve_k + P[m_k, \ve_k] +\dot \ve_k - \theta_k + \thre)
 = \frac{1}{n} \sumkj \dot f_k \dot w_k\,.
\ee
We add \eqref{es5x} to \eqref{ds4}, which yields, by virtue of \eqref{pien},
\begin{eqnarray}\nonumber
\frac{\dd}{\dd t} \left(\frac{1}{n}\sumkj \Big(
\frac 12 \ve_k^2 + V[m_k, \ve_k]+ \thre\ve_k +\frac12 \dot w_k^2 +\theta_k {+ \chi_k}  \Big) + \frac{n}{2} \sumk (\dot w_k - \dot w_{k-1})^2 \right)
\\ = \frac{1}{n} \sumkj (\dot f_k \dot w_k\ + g_k ).\label{es5}
\end{eqnarray}
We estimate the right hand side of \eqref{es5} using the {\em discrete H\"older inequality}
\be{hoe}
\frac{1}{n} \sumk \xi_k\eta_k \le \left(\frac{1}{n} \sumk |\xi_k|^p\right)^{1/p}
\left(\frac{1}{n} \sumk |\eta_k|^{p'}\right)^{1/p'}
\ee
for all vectors $(\xi_1, \dots, \xi_n)$, $(\eta_1, \dots, \eta_n)$, and for $1/p + 1/p' = 1$.}
We have 
$$
\frac{1}{n} \sumkj \dot f_k^2(t) \le \int_0^1 f_t^2(x,t)\dd x\,,
$$
hence, by {\eqref{hoe},} Hypothesis \ref{h1}\,(vii), and Gronwall's lemma,
\be{es5a}
\frac{1}{n}\sumkj \Big( \dot w_k^2 +
\ve_k^2 + \theta_k \Big)(t) + n \sumk (\dot w_k - \dot w_{k-1})^2(t) \le C\,,
\ee
We conclude in particular  that the approximate solutions exist globally and $T_n = T.$


\subsection{The discrete Dafermos estimate}\label{dafe}

 We test \eqref{ds4} by $\theta_k^{-1/3}$ and   we proceed similarly as in {Subsection \ref{dafer}}.
 The integration by parts is replaced by   the elementary inequality
$$
-(x-y)(x^{-1/3} - y^{-1/3}) \ge 3 (x^{1/3} - y^{1/3})^2,
$$
with the choice $x = \theta_k, y = \theta_{k-1}$.  We obtain  for all $t \in [0,T]$ after summing up from $k =1, ..., n-1$ and integrating in time {the following counterpart of \eqref{dafermos}:}
\begin{eqnarray}
\int_0^t \left( 3n \sumk \left(\theta_k^{1/3} - \theta_{k-1}^{1/3}  \right)^2 + \frac{1}{n} \sumkj  \dot\ve_k^2 \; \theta_k^{-1/3}\right)\dd \tau& \le& \int_0^t \frac{1}{n} \sumkj |\dot\ve_k|  \theta_k^{2/3}\dd\tau + \frac{3}{2n}\sumkj \theta_k^{2/3}(t)\nonumber\\ 
&+&C_1 \int_0^t \frac{1}{n} \sum_{k=1}^{n-1} \theta_k^{2/3} \dd \tau.\, 
\label{de1}
\end{eqnarray}
The last two terms on the right hand side are bounded by virtue of \eqref{es5a}.
By {\eqref{hoe}},
$$
\int_0^t \frac{1}{n} \sumkj |\dot\ve_k|  \theta_k^{2/3}\dd\tau \le \left(\int_0^t \frac{1}{n} \sumkj \theta_k^{5/3} \dd\tau\right)^{1/2} \left(\int_0^t \frac{1}{n} \sumkj \theta_k^{-1/3}
\dot\ve_k^2\dd\tau\right)^{1/2},
$$
hence,
\be{de2}
\int_0^t \left(\frac{1}{n} \sumkj  \dot\ve_k^2 \theta_k^{-1/3}
+ 3n \sumk \left(\theta_k^{1/3} - \theta_{k-1}^{1/3}\right)^2\right)\dd \tau
\le C\left(1 + \int_0^t \frac{1}{n} \sumkj \theta_k^{5/3} \dd\tau\right)\,.
\ee
We now apply the inequality \eqref{gagd} as in the formal  case, which implies that
\be{es8}
\int_0^t \left(\frac{1}{n} \sumkj \theta_k^{-1/3} \dot\ve_k^2
+ n \sumk \left(\theta_k^{1/3} - \theta_{k-1}^{1/3}\right)^2\right)\dd \tau  \le C\,.
\ee
Using \eqref{gagd} again 
we obtain
\be{gag3}
\int_0^t \frac{1}{n} \sumkj \theta_k^{8/3}(\tau) \dd\tau \le C\,,
\ee
and, as a consequence of \eqref{es5x},
\be{gag4}
\int_0^t \frac{1}{n} \sumkj \dot\ve_k^2(\tau) \dd\tau \le C\,.
\ee


\subsection{Higher order discrete estimates}\label{high}

We define $\ve_0, \ve_n$ for $k =0, k =n,$ as solutions to the differential equations
\be{bouh}
\left.
\barr{rl}
\ve_0 + P[m_0,\ve_0] + \dot \ve_0 - \theta_0 + \thre &= 0\,,\\
\ve_n + P[m_n,\ve_n] + \dot \ve_n - \theta_n + \thre &= 0\,,
\earr
\right\}
\ee
with initial conditions $\ve_0(0) = \ve_n(0) = 0$. 
{The values of $m_0, m_n$ are chosen as in \eqref{ds5}, where we choose for $\chi$
the natural ``boundary'' conditions compatible with \eqref{bound}, that is,
$$
\chi_0(t) = \chi_1(t),\;\;\; \chi_n(t) = \chi_{n-1}(t).
$$
}
Then \eqref{ds1} holds for all $k=0, \dots, n$, and we have
\be{dife}
\dot u_k - \dot u_{k-1} = \ve_k - \ve_{k-1}
+ P[m_k, \ve_k] - P[m_{k-1},\ve_{k-1}] + \dot \ve_k - \dot \ve_{k-1} - \theta_k + \theta_{k-1}
\ee
for all $k=1, \dots n$. By \eqref{parts} we have
$$
n \sumk (\dot u_k - \dot u_{k-1})^2 = -n \sumkj \dot u_k (\dot u_{k+1} - 2 \dot u_k + \dot u_{k-1}),
$$
hence, by \eqref{ds1}--\eqref{ds2} and \eqref{dife}, 
\begin{align}\no
&n \sumk (\ve_k - \ve_{k-1}
+ P[m_k, \ve_k] - P[m_{k-1},\ve_{k-1}] + \dot \ve_k - \dot \ve_{k-1} - \theta_k + \theta_{k-1})^2\\ \label{ho1}
&\ = \frac{1}{n} \sumkj  (\ddot w_k - \ddot\ve_k - \dot f_k) (\ve_k + P[m_k, \ve_k] + \dot \ve_k - \theta_k + \thre).
\end{align}
This yields, {as a counterpart to \eqref{hor1},
\begin{align}\no
&\frac{\dd}{\dd t} \left(\frac{n}{2}\sumk (\dot w_k - \dot w_{k-1})^2 + \frac{1}{2n} \sumkj \dot\ve_k^2
+ n \sumk (\ve_k - \ve_{k-1})^2\right) + n\sumk \left((\ve_k - \ve_{k-1})^2
+ (\dot \ve_k - \dot \ve_{k-1})^2\right)\\ \no
&\le C n \sumk \left((P[m_k, \ve_k] - P[m_{k-1},\ve_{k-1}])^2 + (\theta_k - \theta_{k-1})^2\right)\\ \no
&
 + \frac{1}{n} \sumkj |\dot f_k|\, |\ve_k + P[m_k, \ve_k] + \dot \ve_k - \theta_k + \thre|\\ \label{ho2}
& + \frac{\dd}{\dd t} \left(\frac{1}{n}\sumkj (\dot w_k - \dot \ve_k) (\ve_k + P[m_k, \ve_k]
 - \theta_k+ \thre)\right)
 + \frac{1}{n}\sumkj |\dot w_k - \dot \ve_k||\dot\ve_k + P[m_k, \ve_k]_t - \dot \theta_k|.
\end{align}
As in \eqref{pt1}, we have $|P[m_k,\ve_k]_t| \le C(1 + \theta_k + \frac{1}{n} \sumkj|\dot\ve_j|),$
and this enables us to}
estimate the terms on the right hand side of \eqref{ho2} as follows:
\begin{equation} \no
\int_0^t \frac{1}{n}\sumkj |\dot w_k - \dot \ve_k||\dot\ve_k + P[m_k, \ve_k]_t 
-\dot\theta_k|\dd\tau
\le C\left(1 + \Big(\!\int_0^t\!\frac{1}{n}\sumkj\dot \theta_k^2\dd\tau\Big)^{1/2}\right),
\end{equation}
\begin{equation}
\int_0^t \frac{1}{n} \sumkj |\dot f_k|\, |\ve_k + P[m_k, \ve_k] + \dot \ve_k - \theta_k + \thre|\dd\tau
\le C,
\end{equation}
\begin{eqnarray}\nonumber
 \frac{1}{n}\sumkj (\dot w_k - \dot \ve_k) (\ve_k + P[m_k, \ve_k]
 - \theta_k+ \thre)& \leq &C\left( 1 +  {\left( \frac{1}{n}\sumkj\dot \ve_k^2\right)}^{1/2} \right)\left(1 + {\left(\frac{1}{n}\sumkj \theta_k^2\right)}^{1/2}   \right)\\
& \leq &  \frac{1}{4n}\sumkj\dot \ve_k^2 +   C {\left(1 + \frac{1}{n} \int_0^t \sumkj \theta_k^2 \dd \tau  \right)}^{1/2},
\end{eqnarray} where we have used  \eqref{es5a}, \eqref{gag3},  \eqref{gag4} and  Hypothesis \ref{h1} (i). {Similarly to \eqref{rhs5},
\eqref{rhs1}, we have by}
Proposition \ref{3p2} and Hypothesis \ref{h1}\,(i) 
\begin{eqnarray}\nonumber
&& \hspace{-14mm} |P[m_k, \ve_k](t) - P[m_{k-1},\ve_{k-1}](t)| 
\\ \nonumber
&=& \left| \int_0^{\infty}\left( \varphi(m_k,r) \stop_r[\ve_k] - \varphi(m_{k-1},r) \stop_r[\ve_{k-1}]   \right) \dd r \right| \\\nonumber
&= &
\left|\int_0^{\infty}\left( \varphi(m_k,r) -\varphi(m_{k-1},r)\right)  \stop_r[\ve_k] \dd r + \int_0^{\infty}\left( \varphi(m_{k-1},r)  (\stop_r[\ve_k]-\stop_r[\ve_{k-1}])   \right) \dd r \right|\\
\label{ho15}
&\le & C\left(  |m_k - m_{k-1}| + \max_{\tau\in [0,t]} |\ve_k(\tau) - \ve_{k-1}(\tau)| \right), 
\end{eqnarray}
where we have by  \eqref{ds5},  Hypothesis \ref{h1}\ (vi) and (iii), \eqref{D-ep}, \eqref{ds6},  \eqref{hoe} and \eqref{gag4} that 
\begin{eqnarray}\nonumber
 |m_k - m_{k-1}|&\leq& C \int_0^t \left( |\dot \chi_k - \dot \chi_{k-1}|
 + \frac{1}{n}\sum_{j=1}^n |\la_{j-k}- \la_{j-k+1}| D_j(t)\right)\\ \nonumber
&\le& C  \int_0^t \left(  |\dot \chi_k - \dot \chi_{k-1}|
+ \frac{1}{n}\sum_{j=1}^n  |\dot \ve_j||\la_{j-k}- \la_{j-k+1}| \right) \dd \tau\\ \no
&\le& C \left(\int_0^t |\dot \chi_k - \dot \chi_{k-1}|  \dd \tau
+ \frac{1}{n^2} \int_0^t \sum_{j=1}^n |\dot \ve_j| \dd \tau\right)\\ \label{ho16}
&\le& C \left( |\chi_k^0 - \chi_{k-1}^0| + \int_0^t |\theta_k - \theta_{k-1}|\dd\tau + \frac{1}{n} \right), 
\end{eqnarray}
 We estimate the initial conditions as in \eqref{eini1}, and
integrating \eqref{ho2} in time we conclude from} the above considerations that
\begin{align}\no
&
 \frac{1}{2n} \sumkj \dot\ve_k^2(t) + n \sumk (\ve_k - \ve_{k-1})^2
+ \int_0^t n \sumk (\dot \ve_k - \dot \ve_{k-1})^2 (\tau)\dd\tau
\\ \no
&\ \le C\Bigg(1 + \frac{1}{n}\sumkj\theta_k^2(t) + \left(\!\int_0^t\!\frac{1}{n}\sumkj\dot \theta_k^2 
\dd\tau\right)^{1/2} + \int_0^t n \sumk (\theta_k - \theta_{k-1})^2(\tau)\dd\tau\\ \no
& + \int_0^t \int_0^\tau n \sumk (\dot \ve_k - \dot \ve_{k-1})^2 (\tau')\dd\tau'\dd\tau\Bigg).
\end{align}
Gronwall's argument {and \eqref{gag3} then yields the following counterpart to \eqref{laste}
\begin{align}\no
&\frac{1}{2n} \sumkj \dot\ve_k^2(t) + n \sumk (\ve_k - \ve_{k-1})^2(t)
+ \int_0^t n \sumk (\dot \ve_k - \dot \ve_{k-1})^2 (\tau)\dd\tau
\\ \label{ho3}
&\ \le C\left(1 + \left(\int_0^t \frac{1}{n}\sumkj\dot \theta_k^2 
\dd\tau\right)^{1/2} + \int_0^t n \sumk (\theta_k - \theta_{k-1})^2(\tau)\dd\tau\right).
\end{align}
We now test \eqref{ds4} by $\theta_k$ and obtain
$$
\frac{\dd}{\dd t} \left(\frac{1}{n}\sumkj\theta_k^2\right) + n \sumk (\theta_k - \theta_{k-1})^2
\le C\left(\frac{1}{n}\sumkj\big( \dot\ve_k^2 + \theta_k +  |\dot\ve_k|(1+\theta_k)\big)\theta_k\right),
$$
where, by \eqref{gag3}
$$
\int_0^t \frac{1}{n}\sumkj \dot\ve_k^2\theta_k\dd\tau \le \left(\int_0^t \frac{1}{n}\sumkj \dot\ve_k^4\dd\tau\right)^{1/2}
\left(\int_0^t \frac{1}{n}\sumkj \theta_k^2\dd\tau\right)^{1/2}
\le C \left(\int_0^t \frac{1}{n}\sumkj \dot\ve_k^4\dd\tau\right)^{1/2},
$$
$$
\int_0^t \frac{1}{n}\sumkj |\dot\ve_k|\theta_k^2\dd\tau \le \left(\int_0^t \frac{1}{n}\sumkj \dot\ve_k^4\dd\tau\right)^{1/4}
\left(\int_0^t \frac{1}{n}\sumkj \theta_k^{8/3}\dd\tau\right)^{3/4}
\le C \left(\int_0^t \frac{1}{n}\sumkj \dot\ve_k^4\dd\tau\right)^{1/4},
$$
hence, {by analogy to \eqref{rht2},}
\be{ho4}
\frac{1}{n}\sumkj\theta_k^2(t) + \int_0^t n \sumk (\theta_k - \theta_{k-1})^2(\tau)\dd\tau
\le C\left(1 + 
\left(\int_0^t \frac{1}{n}\sumkj \dot\ve_k^4\dd\tau\right)^{1/2}\right).
\ee
Finally, we test \eqref{ds4} by $\dot\theta_k$ and obtain from H\"older's inequality that
\begin{equation}
\label{ho5}
\int_0^t \frac{1}{n}\sumkj\dot \theta_k^2(\tau)\dd\tau +  n \sumk (\theta_k - \theta_{k-1})^2(t) \le C\, \left(1 + 
\int_0^t \frac{1}{n}\sumkj \dot\ve_k^4\dd\tau + \int_0^t \frac{1}{n}\sumkj \theta_k^4\dd\tau\right).
\end{equation}
We have using \eqref{gag3}
$$
\int_0^t \frac{1}{n}\sumkj \theta_k^4\dd\tau \le \max_{\tau\in [0,t]}\max_{k=1, \dots, n} \theta_k^{4/3}(\tau)
\int_0^t \frac{1}{n}\sumkj \theta_k^{8/3}\dd\tau \le C\max_{\tau\in [0,t]}\max_{k=1, \dots, n} \theta_k^{4/3}(\tau)
$$
and, by \eqref{gagd} with $q=\infty$, $s=1$, $p=2$, $\vr = 2/3$,
$$
\max_{k=1, \dots, n} \theta_k(\tau) \le C\left(1 + \left(n \sumk (\theta_k - \theta_{k-1})^2(\tau)\right)^{1/3}
\right),
$$
we infer from \eqref{ho5} that
\begin{align}\label{ho6}
\int_0^t \frac{1}{n}\sumkj\dot \theta_k^2(\tau)\dd\tau +  n \sumk (\theta_k - \theta_{k-1})^2(t)
\le C\left(1+ 
\int_0^t \frac{1}{n}\sumkj \dot\ve_k^4\dd\tau\right).
\end{align}
Combining \eqref{ho3} with \eqref{ho4} and \eqref{ho6} yields
\begin{equation}
\label{ho7}\frac{1}{n} \sumkj \dot\ve_k^2(t)
+ \int_0^t n \sumk (\dot \ve_k - \dot \ve_{k-1})^2 (\tau)\dd\tau \le C\left(1+
\left(\int_0^t \frac{1}{n}\sumkj \dot\ve_k^4\dd\tau\right)^{1/2} \right).
\end{equation}
Using the vector notation \eqref{vect}, we have by \eqref{bouh} and \eqref{gag3} that
$$
|\dot\bve(t)|_2^2 = \frac{1}{n} \sumkj \dot\ve_k^2(t) + \frac{1}{n} (\dot\ve_0^2(t)+ \dot\ve_n^2(t))
\le \frac{1}{n} \sumkj \dot\ve_k^2(t) + \frac{C}{n} \left(1 + \sumkj \theta_k^2(t)\right)
\le C + \frac{1}{n} \sumkj \dot\ve_k^2(t)\,,
$$
and we rewrite \eqref{ho7} as ({compare with  \eqref{last}})
\be{ho8}
\max_{\tau\in [0,t]}|\dot\bve(\tau)|_2^2 + \int_0^t |\bfD\dot\bve(\tau)|_2^2\dd\tau
\le C \left(1+  \left(\int_0^t |\dot\bve(\tau)|_4^4\dd\tau\right)^{1/2} \right).
\ee
We estimate the right hand side of \eqref{ho8} using \eqref{gagd} as follows:
$$
|\dot\bve(\tau)|_4 \le C\big(|\dot\bve(\tau)|_2+|\dot\bve(\tau)|_2^{3/4}|\bfD\dot\bve(\tau)|_2^{1/4}\big).
$$
We have $\int_0^t |\dot\bve(\tau)|_2^2\dd\tau \le C$ by virtue of \eqref{gag4}, hence
\begin{align}\no
\int_0^t |\dot\bve(\tau)|_4^4\dd\tau & \le C\max_{\tau \in [0,t]} |\dot\bve(\tau)|_2^2
\left(\int_0^t |\dot\bve(\tau)|_2^2\dd\tau
 + \left(\int_0^t |\dot\bve(\tau)|_2^2\dd\tau\right)^{1/2}
\left(\int_0^t |\bfD\dot\bve(\tau)|_2^2\dd\tau\right)^{1/2}\right)\\
\label{ho9}
& \le C \max_{\tau \in [0,t]}|\dot\bve(\tau)|_2^2 \left(1+\int_0^t |\bfD\dot\bve(\tau)|_2^2\dd\tau\right)^{1/2}.
\end{align}
Combining \eqref{ho8} with \eqref{ho9} yields
\be{ho11}
|\dot\bve(t)|_2^2 + \int_0^t |\bfD\dot\bve(\tau)|_2^2\dd\tau
\le C .
\ee
Therefore there exist a constant $C>0$ such that
\begin{align}\label{ho12}
\frac{1}{n} \sumkj \dot\ve_k^2(t)
+ \int_0^t n \sumk (\dot \ve_k - \dot \ve_{k-1})^2 (\tau)\dd\tau
+ \int_0^t \frac{1}{n}\sumkj (\dot\ve_k^4+\ve_k^4)(\tau)\dd\tau &\le C\,,\\ \label{ho13}
\int_0^t \frac{1}{n}\sumkj(\theta_k^4+ \dot \theta_k^2)(\tau)\dd\tau
+ n \sumk (\theta_k - \theta_{k-1})^2(t) & \le C
\end{align}
for $t \in [0,T]$. By comparison, we also have
\be{ho14a}
\int_0^t n^3 \sumkj (\theta_{k+1} - 2 \theta_k + \theta_{k-1})^2(\tau)\dd\tau \le C\,,
\ee
and similarly for $u_k$.
Finally, we differentiate \eqref{ds1} once in $t$ and test by $\ddot \ve_k$, \eqref{ds2} twice in $t$
and test by $\ddot w_k$, and sum the two equations up. Using \eqref{ho12}--\eqref{ho13} and {treating
the initial conditions as in \eqref{eini2}--\eqref{eini8}, we} get the estimate
\be{ho14}
\frac{1}{n} \sumkj \ddot w_k^2(t) + n \sumk (\ddot w_k - \ddot w_{k-1})^2(t) 
+ \int_0^t \frac{1}{n}\sumkj \ddot \ve_k^2(\tau)\dd\tau \le C\,.
\ee


\section{Proof of Theorem \ref{t2}}\label{proo}

\subsection{Existence}\label{exis}

For a generic sequence $\{\vp_k: k=0,1, \dots, n\}$ we put $\Delta_k\vp = n(\vp_k-\vp_{k-1})$,
and $\Delta_k^2\vp = n^2 (\vp_{k+1}- 2\vp_k+\vp_{k-1})$, and define
piecewise constant, piecewise linear, and piecewise quadratic interpolations
\begin{align}\label{inc}
\bar\vp\on(x) &=
\left\{
\barr{ll}
\vp_k \ & \for x \in \big[\frac{k-1}{n},\frac{k}{n}\big)\,, \ k=1, \dots, n-1\,,\\[2mm]
\vp_{n-1} \ & \for x \in \big[\frac{n-1}{n}, 1\big]\,,
\earr
\right.\\[3mm]
\label{inl}
\hat\vp\on(x) &=
\barr{ll}
\vp_{k-1} + \big(x-\frac{k-1}{n}\big)\Delta_k\vp
\  \for x \in \big[\frac{k-1}{n},\frac{k}{n}\big)\,, & \quad k=1, \dots, n\,,
\earr
\\[3mm]
\label{inq}
\tilde\vp\on(x) &=
\left\{
\barr{ll}
\frac12(\vp_{k-1}+\vp_k) + \big(x-\frac{k-1}{n}\big)\Delta_k\vp
 + \frac{1}{2}\big(x-\frac{k-1}{n}\big)^2 \Delta_k^2\vp
 \ & \for x \in \big[\frac{k-1}{n},\frac{k}{n}\big)\,, \\[2mm]
 &\quad k=1, \dots, n-1\,,\\[2mm]
\frac12(\vp_{n-1}+\vp_n) + \big(x-\frac{n-1}{n}\big)\Delta_n\vp
+ \frac{1}{2}\big(x-\frac{n-1}{n}\big)^2 \Delta_{n-1}^2\vp
\ & \for x \in \big[\frac{n-1}{n}, 1\big]\,.
\earr
\right.
\end{align}

We also define 
\begin{equation}
\la^{(n)}(x,y) = \la_{j-k} \qquad \for  (x,y) \in \left[  \frac{k-1}{n}, \frac{k}{n}\right) \times  \left[  \frac{j-1}{n}, \frac{j}{n}\right).
\end{equation}

For functions $\bar\ve\on$, $\bar\theta\on$, $\bar u\on$, $\bar w\on$, $\hat\ve\on$,
$\tilde\theta\on$, $\tilde u\on$, $\tilde w\on$, we have derived estimates
\eqref{ho12}--\eqref{ho14} that we rewrite in the form
\begin{align}\label{es1}
\big|\bar\ve_t(t)\big|_2^2 + \int_0^t \big|\hat\ve_{xt}(\tau)\big|_2^2\dd\tau
+  \int_0^t \big(\big|\bar\ve_t(\tau)\big|_4^4+ \big|\bar\ve(\tau)\big|_4^4\Big)\dd\tau
& \le C,\\ \label{es2}
\int_0^t \big(\big|\bar\theta_t(\tau)\big|_2^2+ \big|\bar\theta(\tau)\big|_4^4\Big)\dd\tau
+ \big|\hat\theta_x(t)\big|_2^2 & \le C,\\ \label{es3}
\int_0^t \big|\tilde\theta_{xx}(\tau)\big|_2^2\dd\tau & \le C,\\ \label{es4}
\big|\bar w_{tt}(t)\big|_2^2 + \big|\hat w_{xtt}(t)\big|_2^2 + \int_0^t \big|\bar\ve_{tt}(\tau)\big|_2^2\dd\tau & \le C,
\end{align}
and by \eqref{ds2}--\eqref{ds3},
\be{es6}
\int_0^t \big( \big|\tilde w_{xxt}(\tau)\big|_2^2 + \big|\tilde u_{xxt}(\tau)\big|_2^2\big)\dd\tau
\le C.
\ee
System {\eqref{ds1}--\eqref{ds6} has the form
\begin{align}\label{aps1}
\bar u\on_t &= \bar\ve\on + P[\bar m\on, \bar\ve\on]+\bar\ve\on_{t}-(\bar\theta\on -\thre),\\ \label{aps2}
\bar w\on_{t} - \bar\ve\on_{t} &= -\tilde u\on_{xx} + \bar f\on,\\ \label{aps3}
\bar\ve\on & = \tilde w\on_{xx},\\ \label{aps4}
\bar\theta\on_t &= \tilde\theta\on_{xx}   +  \bar m\on_{t} \bar K\on
 + \bar D\on + (\bar\ve\on_t)^2 -\bar\theta\on \bar\ve\on_t  + \bar g\on(\bar\theta\on) - \bar \chi\on_t\,,\\ 
 \label{aph3.4}
\bar\chi\on(x,t) &= \stop_{[0,1]} [\bar\chi\on(0), \bar A\on(x,\cdot)](t),\\ \label{aph3.5}
\bar m\on(x,t) &= \stop_{[0,\infty)} [0, \bar S\on(x,\cdot)](t),\\ \label{aph3.2}
\bar A\on(x,t) &= \int_0^t\frac{1}{\gamma}
\left(\frac{L}{\thre}(\bar\theta\on-\thre) \right)(x,\tau)\,\dd\tau,\\
\label{aph3.8}
\bar S\on(x,t) &= \int_0^t\left( -h(\bar\chi\on_t(x,\tau))
+ \int_0^1 \la\on(x-y) \bar D\on(y,\tau) \dd y \right)(x,\tau)\,\dd\tau,\\ \label{aph3.6}
\bar D\on(x,t) &=
\int_0^{\infty}\varphi(\bar m\on,r)\ \stop_r[\bar\ve\on] (\bar\ve\on -\stop_r[\bar\ve\on] )_t(x,t)\dd r,\\ \label{aph3.7}
\bar K\on(x,t) &= - \frac{1}{2}  \int_0^{\infty} \varphi_m( \bar m\on, r) \stop_r^2[\bar\ve\on] \dd r,
\end{align}
with $\bar\chi\on(0)$ chosen in agreement with \eqref{inid}.} We further have
$$
\int_0^t \big|\hat\ve_{tt}(\tau)\big|_2^2\dd\tau \le \int_0^t \frac{2}{n} \sum_{k=0}^n
\ddot\ve_k^2(\tau)\dd\tau \le \int_0^t \left(\frac{2}{n} \sumkj
\ddot\ve_k^2(\tau) + \frac{2}{n}(\ddot\ve_0^2(\tau) + \ddot\ve_n^2(\tau))\right)\dd\tau. 
$$
By \eqref{bouh}, we have for $k=0$ and $k=n$
$$
\ddot\ve_k^2(\tau) \le C(1 + \theta_k^2(\tau) + \dot\theta_k^2(\tau))\,,
$$
hence
\be{es7}
\int_0^t \big|\hat\ve_{tt}(\tau)\big|_2^2\dd\tau \le C\Big(1 + \int_0^t \frac{2}{n} \sumkj
\big(\ddot\ve_k^2 + \theta_k^2 + \dot\theta_k^2\big)(\tau)\dd\tau\Big) \le C. 
\ee
From \eqref{es1}, \eqref{es7}, and from Sobolev embedding theorems
it follows that there exists $\ve \in W^{1,2}(\Omega_T)$ such that
$\ve_{xt}, \ve_{tt} \in L^2(\Omega_T)$, and a subsequence of $\{\hat\ve\on\}$, still indexed by $n$,
such that
$$
\hat\ve\on \to \ve \ \ \hbox{strongly in } C(\overline\Omega_T)\,,
\quad \hat\ve\on_t \to \ve_t \ \ \hbox{strongly in }  L^p(\Omega_T)
$$
for all $p>1$. Furthermore,
$$
|\bar\ve\on_t - \hat\ve\on_t|^2(x,t) \le |\dot\ve_k - \dot\ve_{k-1}|^2(t)
$$
for $x \in [(k-1)/n, k/n]$, hence
$$
\int_0^t\ipi |\bar\ve\on_t - \hat\ve\on_t|^2(x,\tau)\dd x \dd\tau
\le \int_0^t \frac{1}{n}\sumk (\dot \ve_k - \dot \ve_{k-1})^2 (\tau)\dd\tau \le \frac{C}{n^2}\,,
$$
so that
$$
\bar\ve\on_t \to \ve_t \ \ \hbox{strongly in }  L^2(\Omega_T).
$$
Similarly,
$$
|\bar\ve\on - \hat\ve\on|^2(x,t) \le |\ve_k - \ve_{k-1}|^2(t)
\le \sumk (\ve_k - \ve_{k-1})^2 (t) \le \frac{C}{n},
$$ 
hence
$$
\quad \bar\ve\on \to \ve \ \ \hbox{strongly in } {L^\infty(\Omega_T).}
$$
We check in the same way that there exist $u, w, \theta \in C(\overline\Omega_T)$ such that,
selecting again a subsequence if necessary, 
$$
\tilde w\on_{xxt} \to \ve_t = w_{xxt}\,, \ \tilde u\on_{xx} \to u_{xx} \ \ \hbox{strongly in } L^2(\Omega_T),
$$
$$
\bar\theta_t\on \to \theta_t\,, \ \tilde\theta\on_{xx} \to \theta_{xx} \ \ \hbox{weakly in } L^2(\Omega_T),
\quad \bar\theta\on \to \theta \ \ \hbox{strongly in } L^\infty(\Omega_T).
$$
Finally, {for all $n,l \in \nat$} we have 
\[
|\bar\chi^{(n)}(x,t) -{\bar\chi^{(l)}(x,t)} | \leq 2 \max_{\tau \in [0,t]} |\bar A^{(n)}
- {\bar A^{(l)}}|(x,\tau)
\]
and 
\[
|\bar\chi^{(n)} - \bar\chi^{(l)}|(x,t) \le
\int_0^t |\bar\chi_t^{(n)} - \bar\chi_t^{(l)}|(x,\tau) \dd \tau
\leq C \int_0^t  |\bar \theta^{(n)} - {\bar \theta^{(l)}}|(x,\tau) \dd\tau
+ |\bar\chi^{(n)} - \bar\chi^{(l)}|(x,0). 
\]
It follows that $\bar\chi^{(n)}$ and $\bar\chi_t^{(n)}$ are Cauchy sequences in
$L^\infty(\Omega_T)$ and in {$L^\infty(0,1; L^1(0,T))$,} respectively.
Moreover we have for all $x \in \Omega$ by Proposition \ref{3p2}\,(ii)
that
\begin{eqnarray}\nonumber
\hspace{-26mm}&&\int_0^t \left|  \bar m_t^{(n)} - \bar m_t^{(l)} \right|(x,\tau) \dd \tau  
\leq C \int_0^t |\bar\chi_t^{(n)} -\bar\chi_t^{(l)} |(x,\tau) \dd \tau \\ \label{convm}
\hspace{-16mm} && +\int_0^t \int_0^1 \int_0^{\infty} \left| \la^{(n)}(x,y) \varphi (\bar m^{(n)},r) \delta^{(n)} (y,t,r) -
\la^{(l)}(x,y) \varphi (\bar m^{(l)},r) \delta^{(l)} (y,t,r)\right| \dd r \dd y \dd \tau,
\end{eqnarray}
where we denote
\[
\delta^{(n)} = \delta^{(n)}(y,t,r) = \stop_r[\bar \ve^{(n)}](\bar \ve^{(n)} -\stop_r[\bar \ve^{(n)}])_t (y,t)
{ = r|\play_r[\bar \ve^{(n)}]_t (y,t)|.}
\]
By Proposition \ref{3p2}\,(ii) we have
\[
\int_0^t |\delta^{(n)} - \delta^{(l)}|(y,\tau)\dd\tau \le 
{r \int_0^t |\bar \ve^{(n)}_t - \bar \ve^{(l)}_t|(y,\tau) \dd\tau,}
\]
hence
\[
\int_0^t \int_0^1 \int_0^{\infty} \la^{(n)}(x,y)\varphi(\bar m^{(n)},r) |\delta^{(n)} - \delta^{(l)}|
 \dd r \dd y \dd \tau
\le { C \int_0^t \int_0^1 |\bar \ve^{(n)}_t - \bar \ve^{(l)}_t|(y,\tau) \dd y\dd\tau.}
\]
Similarly,
\begin{align} \no
&\hspace{-14mm}
\int_0^t \int_0^1 \int_0^{\infty}\delta^{(l)} \la^{(n)}(x,y) |\varphi(\bar m^{(n)},r)- \varphi(\bar m^{(l)},r)|
\dd r \dd y \dd \tau\\ \no
& \le C \int_0^t \left(\int_0^1|\bar \ve^{(n)}_t(y,\tau)|\dd y\right)
\max_{x\in\Omega} | m^{(n)}(x,\tau) - m^{(l)}(x,\tau)|\dd \tau.
\end{align}
Finally, we have the pointwise bound
\[
|\la^{(n)}(x,y)- \la^{(l)}(x,y)| \leq \frac{4\Lambda}{\min\{n,l\}.}
\]
We thus have transformed \eqref{convm} into the inequality
\bearr \nonumber 
\max_{x\in\Omega} | m^{(n)} - m^{(l)}|(x,t) &\le&
\max_{x\in\Omega}\int_0^t \left|  \bar m_t^{(n)} - \bar m_t^{(l)} \right|(x,\tau) \dd \tau \\ \label{convmt}
&\le& {q_{nl}} + C \int_0^t \left(\int_0^1|\bar \ve^{(n)}_t(y,\tau)|\dd y\right)
\max_{x\in\Omega} | m^{(n)} - m^{(l)}|(x,\tau)\dd \tau,\qquad
\eearr
{with 
\[
q_{nl} = C\left(\frac{1}{\min\{n,l\}} + |\bar\chi^{(n)}(\cdot,0) - \bar\chi^{(l)}(\cdot,0)|_1
+ \|\bar\theta^{(n)} - \bar\theta^{(l)}\|_\infty
+ \|\bar\ve^{(n)}_t - \bar\ve^{(l)}_t\|_1\right).
\] }
Inequality \eqref{convmt} can} be interpreted as an inequality of the form
\[
q(t) \le {q_{nl}} + \int_0^t s(\tau) q(\tau)\dd\tau,
\]
with $q(t) = \max_{x\in\Omega} | \bar m^{(n)} - \bar m^{(l)}|(x,t)$,
$s(t) = C \int_0^1|\bar \ve^{(n)}_t(y,t)|\dd y$,
$s \in L^1(0,T)$. We obtain using Gronwall's lemma that
$$
q(t) \le {q_{nl}\expe^{\int_0^t s(\tau)\dd\tau} \le C q_{nl},}
$$
so that
$$
\bar m\on \to m \ \ \hbox{strongly in }  L^\infty(\Omega_T),
$$
and, by \eqref{convmt},
$$
\bar m_t\on \to m_t \ \ \hbox{strongly in }  L^\infty(0,1; L^1(0,T)),
$$
This enables us to pass to the limit in
{\eqref{aps1}--\eqref{aph3.7}}
and conclude that $(u,w,\theta,m, \chi)$ is a strong
solution to \eqref{s1}--\eqref{s4} with the regularity indicated in Theorem \ref{t2}
and satisfying the initial conditions \eqref{ini}.
It remains to check that the boundary conditions \eqref{bou} hold.
We have $w_n(t) = 0$, hence
\begin{align}\no
|\tilde w\on(1,t)| &= |2 w_n(t) - \frac32 w_{n-1}(t) + \frac12 w_{n-2}(t)| = 
|w_n(t) - w_{n-1}(t) - \frac12 (w_{n-1}(t) - w_{n-2}(t))|\\ \no
 &\le 2\Big(\sumk |w_k - w_{k-1}|^2(t)\Big)^{1/2} \le \frac{C}{\sqrt{n}},
\end{align}
and similarly for $w(0,t), u(1,t), u(0,t)$. To complete the existence proof, we only have
to verify the homogeneous Neumann boundary condition for $\theta$. In other words, we have to
check that for every $\tilde{\psi}\in C^1(\overline\Omega_T)$ we have
\be{bthe}
\int_0^{T} \ipi (\theta_x\tilde{\psi}_x + \theta_{xx} \tilde{\psi})(x,t)\dd x\dd t = 0\,.
\ee
A straightforward computation yields
\be{bthen}
\int_0^{T} \ipi (\tilde\theta\on_x\tilde{\psi}_x + \tilde\theta\on_{xx} \tilde{\psi})(x,t)\dd x\dd t = 
{- \int_0^T} \tilde{\psi}(1,t) n(\theta_{n-1} - \theta_{n-2})(t)\dd t\,.
\ee
We have
\begin{align}\no
{\int_0^T} n^2(\theta_{n-1} - \theta_{n-2})^2(t)\dd t &=
{\int_0^T} n^2(\theta_n - 2\theta_{n-1} + \theta_{n-2})^2(t)\dd t\\ \no
& \le {\int_0^T} n^2 \sumkj(\theta_{k+1} - 2\theta_{k} + \theta_{k-1})^2(t)\dd t \le \frac{C}{n}\,,
\end{align}
hence, {by virtue of \eqref{ho14a},}
$$
\limn \int_0^{T} \ipi (\tilde\theta\on_x\tilde{\psi}_x + \tilde\theta\on_{xx} \tilde{\psi})(x,t)\dd x\dd t = 0
$$
and \eqref{bthe} follows.


\subsection{Uniqueness}\label{uniq}

Let $(u,w,\theta, \chi, m)$, $(\tilde u, \tilde w, \tilde \theta, \tilde \chi, \tilde m)$
be two solutions of \eqref{s1}--\eqref{bou}, {with the same initial conditions and
the same right hand sides}.
We integrate the difference of \eqref{s3} for $\theta$ and $\tilde\theta$ in time,
{and estimate the terms on the right hand side as follows:
\be{uni1a}
\int_0^t |D[m,w_{xx}] - D[\tilde m,\tilde w_{xx}]|(x,\tau)\dd\tau
\le C\left(\int_0^t (|m - \tilde m| + |w_{xxt} - \tilde w_{xxt}|)(x,\tau)\dd\tau\right),
\ee
where we have used Hypothesis \ref{h1}\,(i) and Proposition \ref{3p2}\,(ii). Furthermore,
\be{uni2a}
\int_0^t |\theta w_{xxt} - \tilde\theta \tilde w_{xxt}|\dd\tau
\le \|\tilde\theta\|_\infty \int_0^t |w_{xxt} - \tilde w_{xxt}|\dd\tau
+ \left(\int_0^t |w_{xxt}|^2\dd\tau\right)^{1/2}
\left(\int_0^t |\theta - \tilde\theta|^2\dd\tau\right)^{1/2}.
\ee
We have by \eqref{rht9} that
\[
\max_{x\in[0,1]}\int_0^t |w_{xxt}|^2(x,\tau)\dd\tau \le C(\|w_{xxt}\|_2^2 + \|w_{xxxt}\|_2^2) \le C,
\]
hence
\be{uni3}
\int_0^t |\theta w_{xxt} - \tilde\theta \tilde w_{xxt}|(x,\tau)\dd\tau
\le C\left(\int_0^t |w_{xxt} - \tilde w_{xxt}|(x,\tau)\dd\tau
+ \left(\int_0^t |\theta - \tilde\theta|^2(x,\tau)\dd\tau\right)^{1/2}\right).
\ee
Similarly,
\be{uni4}
\int_0^t |w_{xxt}^2 - \tilde w_{xxt}^2|(x,\tau)\dd\tau
\le C \left(\int_0^t |w_{xxt} - \tilde w_{xxt}|^2(x,\tau)\dd\tau\right)^{1/2}.
\ee
The fatigue term is estimated as
\[
\int_0^t |m_t \KK[m,w_{xx}] - \tilde m_t \KK[\tilde m,\tilde w_{xx}]|\dd\tau
\le C \int_0^t (|m_t - \tilde m_t|+ |m_t||m - \tilde m|+ |w_{xxt} - \tilde w_{xxt}|)\dd\tau,
\]
where $|m_t(x,t)| \le C$ by virtue of \eqref{rht8}, and
\be{uni5}
\int_0^t |m_t - \tilde m_t|(x,\tau)\dd\tau \le C\int_0^t \left( |\theta - \tilde\theta|(x,\tau)
+ \int_0^1 (|m - \tilde m| + |w_{xxt} - \tilde w_{xxt}|)(y,\tau)\dd y\right)\dd\tau
\ee
by Proposition \ref{3p2}\,(ii). From Gronwall's argument we obtain
{
\be{uni6}
\int_0^t |m_t - \tilde m_t|(x,\tau)\dd\tau \le C\int_0^t \left( |\theta - \tilde\theta|(x,\tau)
+ \int_0^1 (|\theta - \tilde\theta| +
 |w_{xxt} - \tilde w_{xxt}|)(y,\tau)\dd y\right)\dd\tau\,.
\ee}
Finally,
\[
|\chi(x,t) - \tilde\chi(x,t)| \le C \int_0^t |\theta - \tilde\theta|(x,\tau)\dd\tau, \
\int_0^t |g(\theta,x,\tau) - g(\tilde\theta,x,\tau)|\dd\tau
\le C\int_0^t |\theta - \tilde\theta|(x,\tau)\dd\tau.
\]
We now test the resulting inequality} by
$\theta - \tilde\theta$ {and integrate in $x$.} Taking into account the above estimates,
we finally obtain
\begin{align}\no
&\ipi |\theta-\tilde\theta|^2(x,t)\dd x + \frac12\frac{\dd}{\dd t}\ipi\Big(\int_0^t (\theta_x - \tilde\theta_x)
(x,\tau)\dd\tau\Big)^2 \dd x \\ \label{uni1}
&\hspace{20mm}\le C \int_0^t\ipi \big(|w_{xxt} - \tilde w_{xxt}|^2 + |\theta-\tilde\theta|^2\big)(x,\tau)\dd x \dd\tau\,.
\end{align}
In the next step, we test the difference of the time derivatives of \eqref{s2} for $w$ and $\tilde w$
by $w_t-\tilde w_t$, the difference of \eqref{s1} for  $u$ and $\tilde u$ by  $w_{xxt} - \tilde w_{xxt}$,
and sum up. Arguing as above, we obtain
\begin{align}\no
&\frac12\frac{\dd}{\dd t} \ipi\big((w_t- \tilde w_t)^2 + (w_{xt} - \tilde w_{xt})^2\big)(x,t)\dd x
+ \ipi |w_{xxt} - \tilde w_{xxt}|^2(x,t)\dd x\\ \label{uni2}
& \quad \le  C\Big(\ipi |\theta-\tilde\theta|^2 (x,t)\dd x + 
 \int_0^t\ipi |w_{xxt} - \tilde w_{xxt}|^2(x,\tau)\dd x \dd\tau\Big).
\end{align}
{}From \eqref{uni1}--\eqref{uni2} it follows that
\begin{align}\no
&\ipi \big(|w_{xxt} - \tilde w_{xxt}|^2 + |\theta-\tilde\theta|^2 \big)(x,t) \dd x\\ \no
&\hspace{25mm} +
\frac12\frac{\dd}{\dd t}\ipi\Big(\Big(\int_0^t (\theta_x - \tilde\theta_x)
(x,\tau)\dd\tau\Big)^2 + (w_t- \tilde w_t)^2 + (w_{xt} - \tilde w_{xt})^2\Big)(x,t)\dd x \\ \no
&\hspace{20mm}\le C \int_0^t\ipi \big(|w_{xxt} - \tilde w_{xxt}|^2 + |\theta-\tilde\theta|^2\big)(x,\tau)\dd x \dd\tau\,.
\end{align}
Gronwall's argument now yields that $w=\tilde w$, $\theta = \tilde\theta$, and the proof of Theorem
\ref{t2} is complete.

\end{document}